\documentclass[12pt]{article}

\usepackage{amsfonts,amssymb,amsmath,eucal}
\usepackage[dvips]{graphics}

\numberwithin{equation}{section}

\newcommand{\Z}{\mathbb{Z}}
\newcommand{\R}{\mathbb{R}}
\newcommand{\C}{\mathbb{C}}
\newcommand{\Mh}{\widehat{M}}
\newcommand{\bro}{\boldsymbol{\varrho}}
\newcommand{\bsi}{\boldsymbol{\sigma}}

\newtheorem{thm}{Theorem}
\newtheorem{prop}{Proposition}
\newtheorem{lem}{Lemma}
\newcommand{\E}[1]{\left\langle {#1} \right\rangle}
\newcommand{\tr}{\operatorname{tr}}
\newcommand{\Cov}{\operatorname{Cov}}
\newcommand{\Prob}{\operatorname{Prob}}

\newcommand{\Map}{\operatorname{Map}}
\newcommand{\Metr}{\operatorname{Met}}
\newcommand{\Aut}{\operatorname{Aut}}
\newcommand{\Ai}{\operatorname{Ai}}
\newcommand{\map}{\operatorname{map}}
\newcommand{\const}{\operatorname{const}}
\newcommand{\Img}{\operatorname{Im}}

\newcommand{\val}{\operatorname{val}}

\newcommand{\Su}{\mathsf{S}}
\newcommand{\tX}{\widetilde{X}}
\newcommand{\hx}{\widehat{x}}
\newcommand{\hy}{\widehat{y}}

\newcommand{\la}{\lambda}
\newcommand{\al}{\alpha}
\newcommand{\tC}{\widetilde{C}}
\newcommand{\sq}{\square}
\newcommand{\dn}{\delta}
\newcommand{\up}{\delta^*}
\newcommand{\Pl}{\mathfrak{P}}
\newcommand{\eps}{\varepsilon}
\newcommand{\fX}{\mathcal{X}}
\newcommand{\fY}{\mathcal{Y}}

\begin{document}

\title{Random Matrices and Random Permutations}
\author{Andrei Okounkov\thanks{
 Department of Mathematics, University of California at
Berkeley, Evans Hall \#3840, 
Berkeley, CA 94720-3840. E-mail: okounkov@math.berkeley.edu}}
\date{}
\maketitle

\begin{abstract}
We prove the conjecture of Baik, Deift, and 
Johansson which says that with respect to the Plancherel
measure on the set of partitions $\lambda$ of $n$, 
the rows $\la_1,\la_2,\la_3,\dots$ of $\la$ behave, suitably scaled, like
the 1st, 2nd, 3rd, and so on eigenvalues of a Gaussian random 
Hermitian matrix as $n\to\infty$. Our proof is based on an
interplay between maps on surfaces and ramified coverings
of the sphere. We also establish a connection of this problem with
intersection theory on the moduli spaces of curves. 
\end{abstract}

\section{Introduction}

\subsection{Plancherel measures}

The Plancherel measure is probability measure defined on
the set $G^\wedge$ of irreducible
representations $\pi$ of any finite group $G$. Concretely,
the measure of a representation $\pi$ is
$(\dim\pi)^2\big/|G|$. It is called
Plancherel  because the Fourier transform 
$$
L^2(G,\mu_{\mathrm{Haar}})^G \xrightarrow{\quad Fourier\quad } L^2(G^\wedge,
\mu_{\mathrm{Plancherel}})
$$
is an isometry just like in the classical Plancherel theorem. 

In this paper we will be dealing with Plancherel measures for $S(n)$
and their asymptotics as $n\to\infty$.  The set $S(n)^\wedge$ is labeled by
partitions $\la$ of $n$ or, equivalently, by Young diagrams with
$n$ squares. We denote the Plancherel measure by
$$
\Pl_n(\la)=\frac{(\dim\la)^2}{n!}\,,\quad |\la|=n\,, 
$$
and recall that the dimension $\dim\la$ is given by several 
classical formulas such as  the hook formula. 

\subsection{Limit shape} 

\begin{figure}[!hbt]
\centering
\scalebox{.7}{\includegraphics{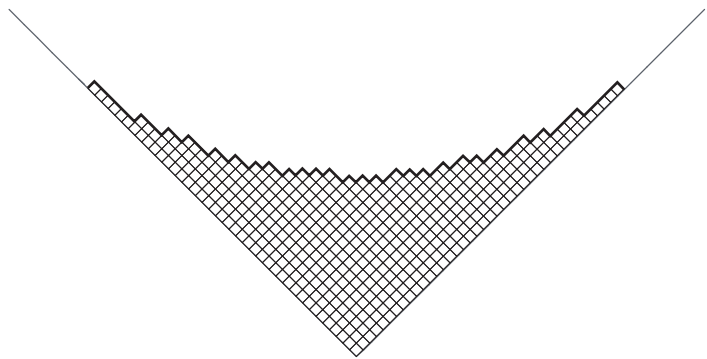}}
\caption{A Young diagram rotated $135^\circ$}
\label{fig1}
\end{figure}

\begin{figure}[!hbt]
\centering
\scalebox{.7}{\includegraphics{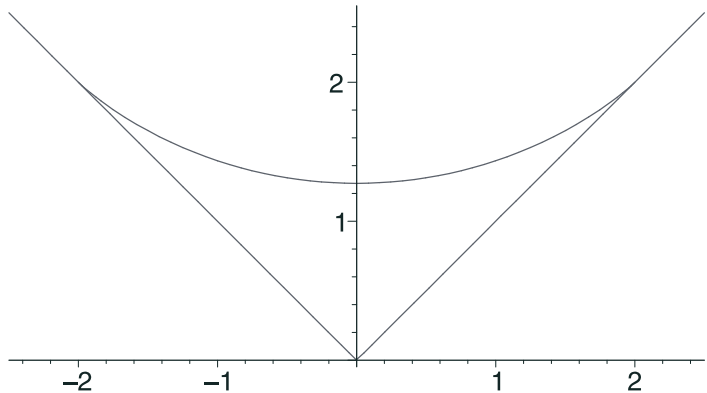}}
\caption{The limit curve $\Omega(x)$}
\label{fig2}
\end{figure}

Logan and Shepp \cite{LS} and, independently, Vershik and Kerov
\cite{VK1} (see also the paper \cite{VK2} which contains
complete proofs of the results announced in \cite{VK1})
discovered the following measure concentration
phenomenon for the Plancherel measures for
$S(n)$ as $n\to\infty$. Take a diagram $\la$, scale it in both directions
by a factor of $n^{-1/2}$ so that to obtain a shape of unit area, 
and rotate it by $135^\circ$ like in Figure \ref{fig1}. 
The boundary of this shape is a polygonal line which is thickened in 
Figure \ref{fig1}. In this way the Plancherel measure $\Pl$ becomes
a measure on the space of continuous functions. It was shown in
\cite{LS,VK1,VK2} that 
as $n\to\infty$, this measure converges to the delta measure
at the following function
\begin{equation}\label{e1}
\Omega(x)=
\begin{cases}
\frac2\pi\left(x \arcsin(x/2)+\sqrt{4-x^2}\right) & |x|\le 2,\\
|x|& |x|\ge 2\,,
\end{cases}
\end{equation}
whose graph is drawn in Figure \ref{fig2}. 

The constant 2 in \eqref{e1} means that the first part of $\la$
should behave like $\sim2\sqrt{n}$ as $n\to\infty$. Indeed,
it was shown in \cite{VK1,VK2} that $\la_1/\sqrt{n}\to 2$ in
probability (in \cite{LS} the inequality $\lim \la_1/\sqrt{n} \ge  2$
was obtained). This constant 2 
corresponds to the constant 2 in the Ulam problem about
the length of the longest increasing subsequence in a
random permutation; it was also obtained by
different means in  \cite{AD,J,S}. About the history
of the Ulam problem see  \cite{AD,BDJ1} and vast literature
cited there. 

\subsection{CLT for limit shape}

The next term in the asymptotic of the Plancherel 
measure was computed by Kerov in \cite{Ke1} who showed that
the Plancherel measure behaves like 
\begin{equation}\label{e2}
\Omega(x)+\frac{U(x)}{n^{1/2}}+ o\left(\frac1{n^{1/2}}\right)\,, 
n\to\infty\,,
\end{equation}
where $U(x)$ is the following Gaussian random process
$$
U(x)=\sum_{k=1}^{\infty} \zeta_k \, \frac{U_k(x)}{\sqrt{k+1}}\,.
$$
Here $U_k(x)$ are the Tchebychef polynomials of the second kind
$$
U_k(2\cos\phi)=\frac{\sin (k+1)\phi}{\sin\phi}
$$
and $\zeta_k$ are independent standard normal variables. Observe
that near the endpoints $x=\pm 2$ the formula \eqref{e2}
becomes inadequate because the series diverges at the endpoints.
For more information about the behavior of the Plancherel typical
partition in the bulk of the limit shape the reader is
referred to the recent paper \cite{BOO}. 

\subsection{Edge of the limit shape and Baik-Deift-Johansson conjecture}

The behavior of the Plancherel measure near the edges of $[-2,2]$ has
been the subject of intense recent studies and numerical experiments,
see \cite{BDJ1} and references therein. 
It has been conjectured by Baik, Deift, and Johansson
that this behavior, suitably scaled, is
identical to the behavior of the eigenvalues of a random Hermitian matrix near
the edge of the Wigner semicircle. More precisely, consider a random  
$n\times n$ matrix
$$
H=
\left(
\begin{array}{ccc}
&\vdots&\\
\hdots&h_{ij}&\hdots\\
&\vdots&
\end{array}
\right)_{1\le i,j\le n}\,, \quad h_{ij}=\overline{h_{ji}}
$$
such that the real and imaginary parts 
$$
h_{ij}=u_{ij}+i v_{ij}
$$
are independent normal variables with mean 0 and variance $1/2$. 
Let 
$$
E_1 \ge E_2 \ge E_3 \ge \dots
$$
be the eigenvalues of $H$. Introduce the variables $y_i$
\begin{equation}\label{e001}
y_i=n^{2/3}\left(\frac{E_i}{2 n^{1/2}}-1\right)\,, \quad i=1,2,\dots \,.
\end{equation}
Then as $n\to\infty$ the $y_i$'s have a limit distribution
which was studied in \cite{F,TW} and other papers.
In particular, the correlation
functions of this random point process have determinantal
form with the Airy kernel, see for example \cite{TW}.  The distributions of 
individual $y_i$'s were obtained by Tracy and Widom in \cite{TW};
they  involve certain solutions of the  Painlev\'e II
equation. 

Similarly, let $\la=(\la_1,\la_2,\dots)$ be a partition and set
(note the difference with \eqref{e001} in the exponent of $n$)
\begin{equation}\label{e011}
x_i=n^{1/3}\left(\frac{\la_i}{2 n^{1/2}}-1\right)\,, \quad i=1,2,\dots\,.
\end{equation}
Baik, Deift, and Johansson conjectured that the limit distribution of the
$x_i$'s exists and coincides with that of the $y_i$'s. They verified
this conjecture for the distribution of $x_1$ and $x_2$ in \cite{BDJ1}
and \cite{BDJ2}, respectively, using very advanced analytic methods.

\subsection{Main result}

The aim of this paper is to give
a direct combinatorial proof of proof of the following result.

Consider the points $x_1,x_2,\dots$ as a random measure on $\R$
with masses 1 placed at the points $x_i$, $i=1,2,\dots$. Consider
its Laplace transform
$$
\hx(\xi)=\sum_{i=1}^\infty  \exp(\xi x_i) \,, \quad \xi>0 \,,
$$
this is a random process on $\R_{>0}$. Define $\hy(\xi)$
similarly. Denote expectation by angle brackets.  

\begin{thm}\label{t1}
In the $n\to\infty$ limit, all mixed moments of the
random variables $\hx(\xi)$ exist and are identical 
to those of $\hy(\xi)$,
that is, 
\begin{equation}\label{e3}
\lim_{n\to\infty} \Big\langle\,\hx(\xi_1)\cdots\hx(\xi_s)\,\Big\rangle = 
\lim_{n\to\infty} \Big\langle\,\hy(\xi_1)\cdots\hy(\xi_s)\,\Big\rangle\,, 
\end{equation}
for any $s=1,2,\dots$ and any numbers $\xi_1,\dots,\xi_s>0$.
\end{thm}

{F}rom Theorem \ref{t1}  one obtains the following result
about the distribution of the individual rows of a
Plancherel typical partition $\la$ 

\begin{thm}\label{t1'} In the $n\to\infty$ limit, the joint
distribution of $x_1,\dots,x_k$ is identical to  the 
joint distribution of $y_1,\dots,y_k$ for any fixed $k$.
\end{thm}

\subsection{Maps on surfaces vs.\  branched coverings} 

In our proof of Theorem \ref{t1} we use the equivalence 
of two points of view on
topological surfaces (or algebraic curves). One way to think
about a surface is to imagine it glued from polygons by
identifying sides of polygons in pairs. Such a representation
is a combinatorial structure called a
{\em map} on a surface. In connection with quantum gravity,
it has been long known that maps are most
intimately related to random matrices, see e.~g.\ \cite{Z}
for an elementary introduction.

Another equally classical way of representing a 
surface is to realize it as a ramified covering
the sphere $S^2$, or in other words, as a Riemann surface
of an algebraic function of one complex variable. It is classically known that
every problem about the combinatorics of covering has a
translation into a problem about permutations which arise as
monodromies around the ramification points.

The two sides of \eqref{e3} have a combinatorial interpretation as 
asymptotics of certain maps and coverings, respectively. We
produce a correspondence between the two enumeration problems
and show that its deviation from being a bijection is negligible
in the $n\to\infty$ limit.

\subsection{Connection to moduli spaces of curves} \label{s14}

The two sides of \eqref{e3} are also very directly 
connected to intersection theory on the moduli
spaces $\overline{\mathcal{M}}_{g,s}$ of genus $g$ 
curves with $s$ marked points.  
Namely,  we show in Section \ref{Brown}
that our enumeration problem for maps (or coverings)
is related to Kontsevich's combinatorial model \cite{K}
for intersection numbers on $\overline{\mathcal{M}}_{g,s}$
by, essentially, a reparametrization. 

This reparametrization involves passage times for the standard 
Brownian motion. As a consequence, our
enumeration asymptotics derived in Theorem \ref{t3},
differs from the unique boxed  formula of \cite{K} 
 by replacing the
Laplace transform variables by their square roots. 

It follows that the limit \eqref{e3} is a close
relative of the so called $s$-point function for
the intersection numbers 
of the $\psi$-classes on $\overline{\mathcal{M}}_{g,s}$.
This can be used to compute the $s$-point
function, see \cite{O4}. 

It is tempting to speculate that both sides of \eqref{e3}
must be certain Riemann integral sums for
the corresponding integrals over $\overline{\mathcal{M}}_{g,s}$
and the only difference between them is that one 
discretizes $\overline{\mathcal{M}}_{g,s}$ using maps
and the other --- using coverings. 

For another application of asymptotics of 
coverings to evaluation of integrals over certain moduli spaces see
\cite{EO}. Another connection between coverings and
moduli spaces was obtained in \cite{ELSV}.

\subsection{Jucys-Murphy elements}

The reader would be hardly surprised to learn that our main
technical tool on the symmetric group side are the 
Jucys--Murphy elements \cite{Ju,M,DG}. In recent years, they have
become all--purpose heavy--duty technical tools in
representation theory of $S(n)$, see for example \cite{B,KO,O1,O2,OV}. 
 The observation that in the $n\to\infty$ the spectral
measures (in the regular representation) of these elements
becomes the Wigner semicircle was made
by P.~Biane in \cite{B2}.

\subsection{Historic remarks} 
The existence of a connection between Plancherel measures
and random matrices has been actively advocated by S.~Kerov,
 see e.g.~\cite{Ke3,Ke4,Ke5}. The simplest evidence of
such a connection is the fact that the so called  transition
distribution  for the limit shape $\Omega$ coincides with the
Wigner semicircle. Random matrices also enter the 
representation theory of symmetric groups  via the 
free probability theory.  For a detailed discussion of the interplay
between symmetric groups and free probability see the
paper  \cite{B}  by P.~Biane. Our results explain, 
at least to some extend, this connection.  

\subsection{Further development} 

An analytic proof of the Baik--Deift--Johansson conjecture was
found subsequently in \cite{BOO} and, independently, in \cite{J2}.
This approach is based on an exact formula for 
the so-called correlation
functions of the poissonized Plancherel measure.
Same formula allows to analyze the local structure of
a Plancherel typical partition in the bulk of the limit
shape, see \cite{BOO}.  This exact formula
for correlation functions is a special case of the result
obtained in \cite{BoO}. The results of 
\cite{BoO} were considerably generalized in \cite{O3}.

\subsection{Acknowledgements}
The author would like to thank A.~Eskin for numerous discussions and
P.~Biane for explaining the connection to Brownian motion. 
The author was supported by  NSF for under
grant DMS--9801466.  

\section{Random Matrices}   

\subsection{Maps on surfaces and Random Matrices}

\subsubsection{}

The relation between maps on surfaces and random matrices via
the Wick formula is
well known. Classical examples of exploiting this relation
are, for example, the papers \cite{HZ,K}. A very accessible
introduction can be found, for example, in \cite{Z}. See also,
for example, \cite{FGZ} for a physical survey. We briefly
 recall some basic things
in order to facilitate the comparison with 
the enumeration of coverings. 

\subsubsection{}

Consider a random Hermitian matrix  
$$
H=
\left(
\begin{array}{ccc}
&\vdots&\\
\hdots&h_{ij}&\hdots\\
&\vdots&
\end{array}
\right)_{1\le i,j\le n}\,, \quad h_{ij}=\overline{h_{ji}}
$$
such that the real and imaginary parts 
$$
h_{ij}=u_{ij}+i v_{ij}
$$
are independent normal variables with mean 0 and variance $1/2$.
We will be interested in the asymptotics of  
\begin{equation}\label{e21}
\frac1
{2^{|k|}n^{|k|/2}}\E{\prod_{j=1}^s \tr H^{k_i}}\,, \quad k_i\sim \xi_j n^{2/3}
\end{equation}
as $n\to\infty$ and some fixed $\xi_1,\dots,\xi_s>0$.
Here $|k|=\sum_i k_i$.  
Similar averages were considered by many authors, see especially the 
recent paper \cite{Sosh} and references therein. 
Remark that by \eqref{e001} we have 
$$
\left(\frac{E_i}{2 n^{1/2}}\right)^{\xi_j n^{2/3}} \to \exp(\xi_j y_i)\,,
\quad n\to\infty\,,
$$
and that it is clear that only the eigenvalues near the edges of 
the Wigner's semicircle contribute to the asymptotics of 
\eqref{e21}.

\subsubsection{}\label{rel} 

By symmetry, the expectation  \eqref{e21}
vanishes if $|k|$ is odd.   If $|k|$ is even then 
then it is a sum of $2^{s-1}$ terms coming from
various combinations of the maximal and minimal
eigenvalues of $H$. In what follows,
we will always assume that $|k|$ is even. In this
case, there are $2^{s-1}$ possible choices of parity of each
individual $k_i$ and it is easy to see that by taking
a suitable linear combination we can single out
the contribution of only maximal eigenvalues.
Therefore, instead of working with expectations like 
\begin{equation}\label{e002}
\Big\langle\,\hy(\xi_1)\cdots\hy(\xi_s)\,\Big\rangle
\end{equation}
we can work with expectations \eqref{e21} which is more
convenient.  

\subsubsection{Correlation functions} 

Let $\bro(x_1,\dots,x_k)$ denote the $k$-point correlation functions for
the scaled eigenvalues $\dfrac{E_i}{2\sqrt{n}}$ of $H$. 
The expectations \eqref{e21} are closely related to
these correlation functions. 
Let 
$$
\bsi = \sum_i \delta_{E_i/2\sqrt n} \,, 
$$
be the scaled spectral measure of $H$. It is a random measure on $\R$. 
We have
\begin{equation}\label{Hsi}
\frac1
{2^{|k|}n^{|k|/2}}\E{\prod_{j=1}^s \tr H^{k_i}} = \int_{\R^s} u_1^{k_1} \dots u_s^{k_s} \,
\E{\bsi^{\times s}}(du) 
\end{equation}
where $\E{\bsi^{\times s}}$ is the following nonrandom measure on $\R^s$
$$
\E{\bsi^{\times s}}(A_1\times \cdots \times A_s) = \E{\prod_{i=1}^s \bsi(A_i)}\,.
$$

Let $\Pi_s$ be the set of all partitions of the 
set $\{1,\dots,s\}$ into disjoint union of subsets. For any $\al\in\Pi_s$,
denote by $\ell(\al)$ the number of parts in $\al$. For example, 
$$
\al=\{\{1\},\{2,3\}\} \in \Pi_3 \,, \quad \ell(\al)=2\,.
$$
For any $k\in\R^s$ and $\al\in\Pi_s$ denote by $k_\al\in \R^{\ell(\al)}$
the vector with coordinates $\sum_{j\in\al_i} k_j$, where $\al_i$
are the parts of $\al$. We have 
\begin{equation}\label{sigma_rho}
\int_{\R^s} u^k \,
\E{\bsi^{\times s}}(du)  = 
\sum_{\al\in\Pi_s} \int_{\R^{\ell(\al)}}   u^{k_\al} \,
\bro(u_1,\dots,u_{\ell(\al)}) \, du \,.
\end{equation}
For example, for $s=2$ we have 
\begin{equation*}
\int_{\R^2} u_1^{k_1} u_2^{k_2} \,
\E{\bsi^{\times s}}(du)  = \int_{\R^2} u_1^{k_1} u_2^{k_2} \,
\bro(u_1,u_2) \, du_1\, du_2 + \int_{\R^1} u_1^{k_1+k_2} \, \bro(u_1) \, du_1 \,.
\end{equation*}

\subsubsection{Wick formula} 
{F}rom the Wick formula one obtains
\begin{equation}\label{e22}
\frac1
{2^{|k|}n^{|k|/2}}\E{\prod_{j=1}^s \tr H^{k_i}} = \frac1{2^{|k|}} \sum_{\Su} 
n^{\chi(\Su)-s} \, \left|\Map_\Su(k_1,\dots,k_s)\right| \, ,
\end{equation}
where the sum is over all homeomorphism classes of
surfaces $\Su$, not necessarily
connected, $\chi(\Su)$ is the Euler characteristic of the surface $S$, and
$\Map_\Su(k_1,\dots,k_s)$ is the set of solutions to the
the following combinatorial problem. 

\begin{figure}[!hbt]
\centering
\scalebox{.7}{\includegraphics{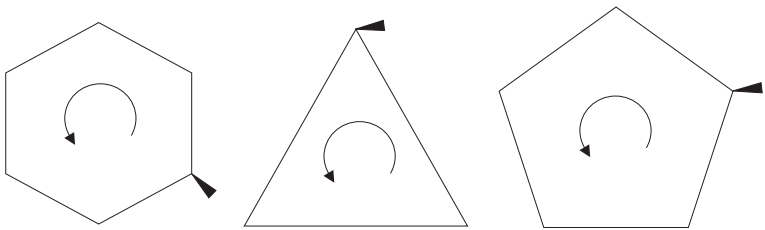}}
\caption{Polygons with orientation and marked vertices}
\label{fig2b}
\end{figure}

Take $s$ polygons:
a $k_1$-gon, a $k_2$-gon, and so on. 
Fix their orientations and a mark a vertex on each
as in Figure \ref{fig2b} (we mark a vertex to
distinguish a $k$-gon from its $(k-1)$ rotations).
Now consider all possible ways to glue 
their sides in pairs in a way consistent with orientation.
The set $\Map_\Su$ consists of all glueings which
produce a surface
homeomorphic to $\Su$.

Note that definition of a map is different 
from the more common one which does not require the choice
of marked vertices. Our definition is more
convenient for our purposes.

\subsubsection{} 

We are interested in the limit of \eqref{e22} as $n$ and the $k_i$'s go
to infinity in such a way that $k_i\propto n^{2/3}$. 
This limit can be determined by taking the term-wise asymptotics in the 
right-hand side of \eqref{e22}. 

For $s=1$ the necessary bounds for the 
 validity of term-wise limits will be obtained from 
explicit formulas in Section \ref{s2ex}. For $s>1$, we then can use
the equation \eqref{sigma_rho} and 
the fact that the correlation functions $\bro$
have a determinantal form and hence
\begin{equation}\label{rhorho}
\bro(u_1,\dots,u_k) \le \prod_{i=1}^k  \bro(u_i) \,.
\end{equation}
Here we use the well-known fact that  the determinant of a positive-definite
matrix is at most the product of its diagonal entries (equivalently,
the volume of a parallelepiped is at most the product of its side lengths).
Consequently, we obtain the estimate \eqref{HZest} which
says that left-hand side of \eqref{e22} is  bounded by
some function of  the variables $\xi_i=k_i\, n^{-2/3}$. 

Since all terms in the right-hand side of
\eqref{e22} are nonnegative, replacing $n$ by a multiple of $n$ we see
that the terms in  in the right-hand side of \eqref{e22} decay 
faster than any exponential and, in particular, taking term-wise limit is
justified.

\subsubsection{} 

Our present goal is to understand the asymptotics of 
$\left|\Map_\Su(k_1,\dots,k_s)\right|$ as the $k_i$'s go to the
infinity. It is clear that it suffices to consider this
asymptotics for connected surfaces $\Su$ only. 
If $\Su$ is a connected surface of genus $g$ we  write
$\Map_g$ instead of $\Map_\Su$.

Below we will describe a function $\map_g(\xi_1,\dots,\xi_s)$ such that
\begin{equation}\label{map}
2^{-|k|} \, \left|
\Map_g(k_1,\dots,k_s)
\right| \sim \map_g(\xi)\, t^{3g-3+3s/2}  \,, \quad \xi_i=k_i/t\,,
\end{equation}
as $t\to\infty$ provided $|k|$ is even. Recall that if $|k|$ is odd
then $\Map_g$ is empty. One extends $\map_\Su$ to disconnected 
surfaces multiplicatively. It is clear from \eqref{map} that 
$\map_\Su$ is homogeneous of total
degree $\frac32(s-\chi(\Su))$ and also positive for
positive values of $\xi$.

It follows that if all
$k_i$'s are even then 
\begin{equation}\label{e23}
\frac1
{2^{|k|}n^{|k|/2}}\E{\prod_{j=1}^s \tr H^{k_i}} \to \sum_{\Su} 
\map_\Su(\xi_1,\dots,\xi_s)\,, \quad k_i\sim \xi_i \, n^{2/3}  \,,
\end{equation}
and if some of the $k_i$'s are odd then  in the right-hand
side of the above formula those terms  that violate the parity conditions 
should be omitted. 

\subsubsection{}
The function $\map_g(\xi)$ can be expressed in terms of the
Laplace transform of the corresponding limits of the correlation
functions $\bro$. 

It is known (see \cite{TW} and note that our $y_i$'s
differ from the centered and scaled eigenvalues which are used in \cite{TW}
by a factor of 2) that 
\begin{equation}\label{rholim}
n^{-2s/3} \, \bro\left(1+\frac{y_1}{n^{2/3}},\dots,1+\frac{y_s}{n^{2/3}}\right)
\to \rho(y_1,\dots,y_s)\,, \quad n\to\infty\,,
\end{equation}
where $\rho$ is given by a determinant 
$$
\rho(y_1,\dots,y_s) = \det \big( K(y_i,y_j)\big) 
$$
with the Airy kernel
$$
K(x,y)= \frac{\Ai(2x) \, \Ai'(2y) - \Ai'(2x) \Ai(2y)}{x-y} \,. 
$$ 
Here $\Ai(x)$ is the classical Airy function. By the l'H\^ospital's rule
and the equation $\Ai''(x)=x \Ai(x)$, we have 
\begin{equation}\label{rho1}
K(x,x)=2 \Ai'(2x)^2 - 4 x \Ai(2x)^2 \,.
\end{equation}

Denote by $R(\xi)$ the Laplace transform
$$
R(\xi_1,\dots,\xi_s)= \int_{\R^s} e^{(\xi,y)} \, \rho(y_1,\dots,y_s) \, dy \,,
$$
which converges for all $\xi\in \R_{>0}^s$. Introduce the function
\begin{equation}\label{defH}
H(\xi_1,\dots,\xi_s) = \sum_{\al\in\Pi_s} R(\xi_\al)\,,
\end{equation}
where, we recall, $\xi_\al$ is the  $\ell(\al)$-dimensional
vector formed by sums of $\xi_i$ over $i$ in blocks of $\al$. 
For example,
\begin{multline*}
H(\xi_1,\xi_2,\xi_3) = R(\xi_1,\xi_2,\xi_3)+ R(\xi_1+\xi_2,\xi_3)+ \\
R(\xi_1+\xi_3,\xi_2) + R(\xi_2+\xi_3,\xi_1)+ R(\xi_1+\xi_2+\xi_3)\,.
\end{multline*}
Finally, set
\begin{equation}\label{defG}
G(\xi_1,\dots,\xi_s) = \sum_{S\subset\{1,\dots,s\}} H(\xi_i)_{i\in S} \,
H(\xi_i)_{i\notin S} \,,
\end{equation}
where the summation is over all subsets $S$ and $H(\xi_i)_{i\in S}$ denotes
the function $H$ in variables $\xi_i$, $i\in S$. For example
$$
G(\xi_1,\xi_2) = 2 H(\xi_1,\xi_2) + 2 H(\xi_1) H(\xi_2) \,.
$$ 

We have from \eqref{e22}, \eqref{Hsi}, \eqref{rholim}, and \eqref{HZest}
$$
G(\xi_1,\dots,\xi_s) =  \sum_{\Su} 
\map_\Su(\xi_1,\dots,\xi_s) \,.
$$
The summation over partitions $\al$ in \eqref{defH} corresponds
to the summation over partitions in \eqref{sigma_rho}. 
The summation over subsets $S$ in \eqref{defG} correspond to the
fact that both ends of the spectrum contribute to the
asymptotics and that the correlations between eigenvalues near
opposite ends of the spectrum disappear in the $n\to\infty$ limit. 

\subsection{Example: maps on the sphere with 1 cell}\label{scat}

\subsubsection{}

\begin{figure}[!hbt]
\centering
\scalebox{.7}{\includegraphics{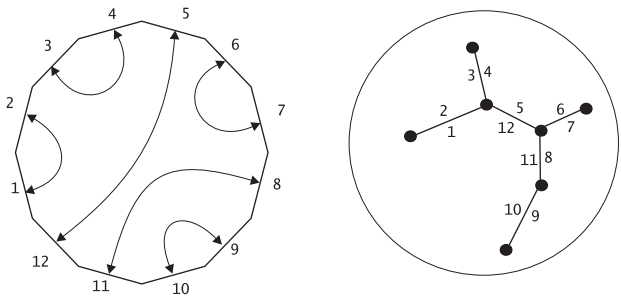}}
\caption{A map on the sphere}
\label{fig3}
\end{figure}

As the simplest example, consider the case $g=0$ and $s=1$, that is, we 
want to glue a sphere from a $k$-gon. One can see that
one obtains a sphere if and only if lines connecting the 
identified sides do not intersect (in which case the
boundary of the polygon becomes a tree in the sphere), see
Figure \ref{fig3}.

The number of such {\em noncrossing} pairings is the Catalan number 
$$
\left|\Map_0(2k)\right|= C_k
=\frac1{k+1}\binom{2k}{k} \,.
$$
The Stirling formula gives
\begin{equation}\label{e003}
\map_0(\xi)=\frac1{\sqrt{\pi}}\left(\frac{\xi}2\right)^{-3/2} \,.
\end{equation}
We have
\begin{equation}\label{ymap}
\E{\hy(\xi)} =\int_{\R^1} e^{\xi y} \,\rho(y) \, dy= \frac12 \sum_{g=0}^\infty \map_g(\xi) \,,
\end{equation}
where 
$$
\rho(y)=K(y,y)
$$ 
is the 1-point correlation functions for the $y_i$'s, that is,
 $\rho(y) \, dy$ is the probability to find one of the $y_i$'s in the interval
$[y,y+dy]$. A formula for this $1$-point function is given in
\eqref{rho1}.

Assuming that we already know that the degree 
of $\map_g(\xi)$ is positive for $g>0$, we conclude that
\begin{equation}\label{as1} 
\int e^{\xi y} \,\rho(y) \, dy \sim \sqrt{\frac{2}{\pi}} \, \frac1{\xi^{3/2}} \,, \quad \xi\to+0\,.
\end{equation}
\begin{figure}[!hbt]
\centering
\scalebox{.7}{\includegraphics{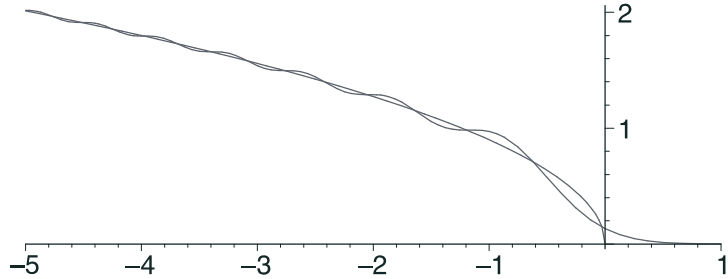}}
\caption{Density $\rho$ of the $y_i$'s versus $2^{3/2}\sqrt{-x}/\pi$}
\label{fig3b}
\end{figure}
This asymptotics reflects the $y\to -\infty$ asymptotics of $\rho(y)$ which is
known to be, see Figure \ref{fig3b},
$$
\rho(y)\sim \frac{2^{3/2}}{\pi}\sqrt{-y}\,, \quad y\to-\infty\,.
$$
This is in agreement with \eqref{as1}.

\subsubsection{}

Below we will need the following elementary lemma about Catalan numbers
\begin{lem}\label{estCat} We have 
$$
 C_k < \frac1{\sqrt{\pi}} \frac{2^{2k}}{k^{3/2}}\,, \quad 
k=1,2,\dots\,.
$$
That is, Catalan numbers $C_k$ are less than their $k\to\infty$ asymptotics. 
\end{lem}

To verify this, set
$$
\tC_k = C_k \, \frac{k^{3/2}}{2^{2k}}  \to  \frac1{\sqrt{\pi}}\,, \quad k\to\infty \,.
$$
We claim that the sequence
$\tC_k$ is strictly increasing, which is equivalent to 
$$
\frac{\tC_{k+1}}{\tC_k} = \frac{(k+1)^{3/2}}{k^{3/2}}\, \frac{k+\frac12}{k+2} > 1 
$$
for $k\ge 1$. Indeed this ratio tends to $1$ as $k\to\infty$ and its
derivative in $k$ is negative: 
$$
\left(\frac{(k+1)^{3/2}}{k^{3/2}}\, \frac{k+\frac12}{k+1}\right)' = 
-\frac34 \,
\frac{(k+1)^{1/2}}{k^{5/2}} \, \frac{3k+2}{(k+2)^2} \,.
$$
Thus, $\tC_k$ is strictly less than the limit $\pi^{-1/2}$ as was to
be shown. 

\subsubsection{}

Another special example to consider is the case $s=2$, $g=0$.
These are the two cases not covered by the general construction 
explained in the next subsection. 

\subsection{Counting maps}

\subsubsection{The contraction $\Phi$} 

We will now count the the maps in all cases except
$s=1,2$, $g=0$. In fact, in order to establish connection 
with random
permutations it is not necessary to actually compute the asymptotics
explicitly.
It suffices to establish just the general pattern of the
combinatorial enumeration which occurs. Nonetheless, we 
do the computations because in the end we will be rewarded
with a connection to the moduli spaces of curves. 

In order to count the maps, we will construct 
a function $\Phi$ from the set of
maps to a simpler set such that the level sets of $\Phi$
are easy to understand. This is like computing the
volume by integrating first along the fibers of a projection
and then over the base. More concretely, the target set of
our function $\Phi$ will be set of pairs
$$
\Phi: \Map_g(k_1,\dots,k_s) \to \{(\Gamma,\ell)\}\,, \quad \Gamma\in
\Gamma^{\ge 3}_{g,s}\,, 
$$
where $\Gamma^{\ge 3}_{g,s}$ denotes \emph{ribbon graphs} of genus $g$ with
$s$ marked cells and vertices of valence $\ge 3$, and 
$\ell$ is a \emph{metric} on the boundary $\partial\Gamma$ of $\Gamma$.

\subsubsection{} 

Recall that, by definition,   a ribbon
graph  is the following object. It is a union of vertices (which are
small disks or polygons; we shall paint them grey in the 
figures) which are connected by ribbons (edges). The boundary
$\partial\Gamma$ of a ribbon graph $\Gamma$ 
is an ordinary graph whose edges are the 
borders of the ribbons. Let $s$ be the number of
connected components of $\partial\Gamma$. Filling each component
of $\partial\Gamma$ with a disk produces a closed surface. The genus of $\Gamma$ is,
by definition, the genus of that surface. 
The components of $\partial\Gamma$ (or the disks filling them)
are called the {\em cells}  of $\Gamma$. We shall consider 
ribbon graphs with $s$ cells and the cells will be marked by
the numbers $\{1,2,\dots,s\}$. 

\subsubsection{} 

Given a map on a surface $\Su$, consider the the graph on $\Su$
formed by vertices and edges of the original polygons. A small 
neighborhood of this graph is a ribbon graph. 
The numbering of the cells comes from the numbering of the
polygons of the map. 

\begin{figure}[!hbt]
\centering
\scalebox{.7}{\includegraphics{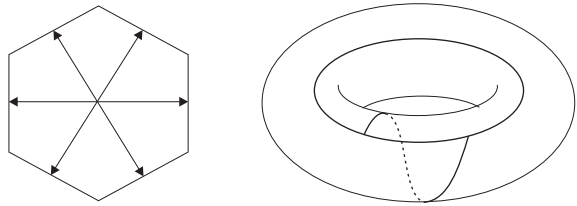}}
\caption{A map on the torus}
\label{fig4}
\end{figure}

\begin{figure}[!hbt]
\centering
\scalebox{.7}{\includegraphics{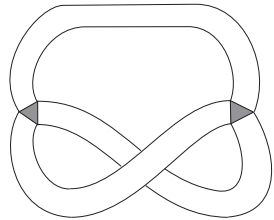}}
\caption{The corresponding ribbon graph}
\label{fig5}
\end{figure}

For example, consider the  map on the 
torus which is drawn in Figure \ref{fig4}. 
The corresponding  ribbon graph is displayed in Figure \ref{fig5}.
There is only one cell in this example.  

We equip the boundary of this graph with the metric $\ell$ in which all
edges have unit length. Thus, we associated
 to any map a pair $(\Gamma_0,\ell)$,
where $\Gamma_0$ is ribbon graph of genus $g$ with $s$ marked cells 
and $\ell$ is a metric on $\partial\Gamma_0$. This pair is the first step
in the construction of $\Phi$. 

\subsubsection{}

The second step on the construction of $\Phi$ is the 
elimination of all  vertices of
valence $\le 2$ from $\Gamma_0$. This goes as follows. 

\begin{figure}[!hbt]
\centering
\scalebox{.7}{\includegraphics{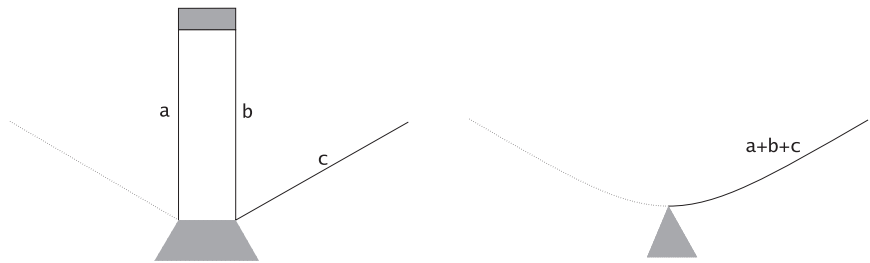}}
\caption{Collapsing of univalent vertices}
\label{fig6}
\end{figure}

\begin{figure}[!hbt]
\centering
\scalebox{.7}{\includegraphics{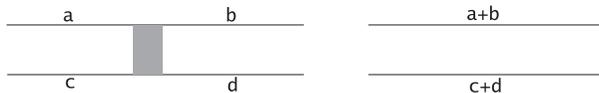}}
\caption{Elimination of  2-valent vertices}
\label{fig7}
\end{figure}

First, we collapse the 
 univalent vertices as in Figure
\ref{fig6}. The numbers in that picture illustrate what we do 
with the metric $\ell$. Namely, we increase the length of 
the adjacent (with respect to the orientation) part of the boundary by the total perimeter of
the disappearing edge. Note that this operation preserves 
perimeters of cells.  
After that, the vertices of valence 2 are eliminated
as in Figure \ref{fig7}. Again, the perimeters of cell
are preserved by this operation. 
In the end, we get a ribbon graph $\Gamma\in\Gamma^{\ge 3}_{g,s}$,
provided we are not in the exceptional cases $g=0$, $s=1,2$ 
in which  we get a point and circle, respectively. 
We also get a metric $\ell$ on $\partial\Gamma$. 
By definition, this pair $(\Gamma,\ell)$ is where $\Psi$
takes our original map. 

\subsubsection{} 

Note that, by construction,
the perimeters of the cells of $\Gamma$ are  equal to  $k_1,k_2,\dots,k_s$.
Also,  the computation of the Euler characteristic gives 
\begin{equation}\label{e007}
\sum_{v\in v(\Gamma)} (\val(v)-2) = 4g-4+2s \,,
\end{equation}
where $v(\Gamma)$ is the set of vertices of $\Gamma$.

All this is, of course, very similar to the stratification of
the moduli space of curves of genus $g$ with $s$ marked points
by means of Strebel differentials, see e.~g.\ \cite{K}. 

\subsection{The level sets of the contraction $\Phi$}\label{sfor}

\subsubsection{}

We now want to compute how many maps $\Phi$ takes to a given pair
$(\Gamma,\ell)$. First, look at a single  edge of $\Gamma$ 
let $p$ and $q$ the lengths of its two boundaries in metric $\ell$. 
We want to compute
how many different configurations produce this data after
the elimination of vertices of valence $\le 2$. 

\begin{figure}[!hbt]
\centering
\scalebox{.7}{\includegraphics{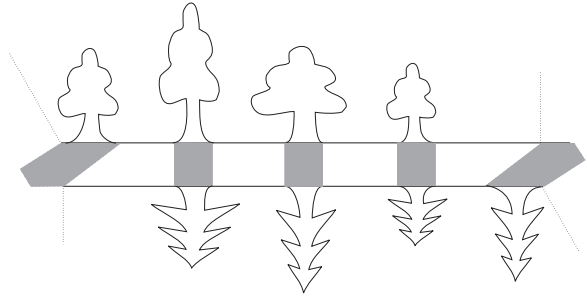}}
\caption{A ribbon graph which collapses to an edge}
\label{fig8}
\end{figure}

This means that we must compute the number
of ribbon graphs of the form shown in Figure \ref{fig8} 
with the length of the upper boundary and lower boundary
being $p$ and $q$, respectively.
The trees in Figure \ref{fig8}
stand for (possibly empty) ribbon graphs which disappear after
collapsing all univalent vertices.  It implies that they are \emph{trees} in
the usual sense of graph theory. 

\begin{figure}[!hbt]
\centering
\scalebox{.7}{\includegraphics{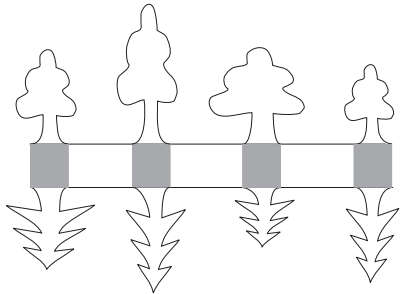}}
\caption{A shift of 
Figure \ref{fig8}}
\label{fig9}
\end{figure}

Remark that a tree is not allowed at one of the ends of both
upper and lower boundary. This corresponds to our convention
(see Figure \ref{fig6}) on where we transfer the length of a collapsing edge.  
However, a simple shift as in Figure \ref{fig9}
which reduces the length of both boundaries by 1, takes 
care of this inconvenience. Now it clear that to obtain
a ribbon graph like in Figure \ref{fig9} one just takes any map 
from $\Map_0(p+q-2)$ and calls the first $(p-1)$ sides
the upper boundary, and rest --- the lower boundary. Therefore,
we get a Catalan number  provided $p+q$ is even
(and $0$ otherwise). This means that 
there are
\begin{equation}\label{e005}
\sim \frac1{\sqrt{2\pi}}\frac{2^{p+q}}{(p+q)^{3/2}}\,,\quad p+q\to\infty\,,
\end{equation}
ribbon graphs which collapse to an edge with length of
the upper and lower boundary equal to $p$ and $q$,
respectively. Moreover, by Lemma \ref{estCat} the actual number of
such maps is always less than \eqref{e005}. 

\subsubsection{} 

Now consider all edges of $\Gamma$. It is clear that we can apply
the above construction to every edge of $\Gamma$ independently and
the only situation in which we get identical maps is when the maps differ by
an automorphism of $\Gamma$. Here by an automorphism we mean automorphisms
of the whole structure of a ribbon graph with marked cells; in particular,
the automorphisms must preserve cells. 

Also recall that we consider maps with marked vertices. Since marks
can be chosen arbitrarily on the boundary of each cell, we have 
\begin{equation}\label{e006}
\left|\Phi^{-1}((\Gamma,\ell))\right| \sim 
\frac{\prod k_i}{|\Aut(\Gamma)|} \, 
\frac{2^{|k|}}{(2\pi)^{|e(\Gamma)|/2}} \, 
\prod_{e\in e(\Gamma)} (\ell_{1,e}+\ell_{2,e})^{-3/2} \,.
\end{equation}
where the $k_i$ are the perimeters of the cells of $\Gamma$
and their product is the number of choices for the marked vertices. 
Here $\ell_{1,e}$ and $\ell_{2,e}$ are the lengths of the two sides
of the edge $e\in e(\Gamma)$ in the metric $\ell$. 

Again, by Lemma \ref{estCat} the right-hand side of \eqref{e006} is
both asymptotics and an estimate from above for the left-hand side. 

\subsection{The asymptotics}

\subsubsection{} 

Now we want to sum \eqref{e006} over the  metrics $\ell$. This means
summation over points in $\R^{2|e(\Gamma)|}$ satisfying the
following properties
\begin{itemize}
\item the values of $\ell$ are integers and for any  $e\in e(\Gamma)$
the sum $\ell_{1,e}+\ell_{2,e}$ is an even integer,
\item the lengths of the edges are nonnegative and 
the perimeters of the $s$ cells are equal to
$k_1,\dots,k_s$, respectively. 
\end{itemize}
It is clear that the first condition defines a sublattice $\Lambda$ of 
index $\left|e(\Gamma)\right|$ in $\Z^{2|e(\Gamma)|}$. 
The second condition defines a convex polytope which we 
denote by $\Metr_\Gamma(k)$. The dimension of this polytope is 
\begin{equation}\label{dimM}
\dim\Metr_\Gamma(k)=2|e(\Gamma)|-s \,.
\end{equation}

By definition, let $\Map_{g,\Gamma}(k)$ denote those maps in $\Map_g(k)$  
which correspond to a given graph $\Gamma$ under $\Phi$. It follows that
\begin{equation}\label{sumM}
\left|\Map_{g,\Gamma}(k)\right| = \sum_{\ell\in \Metr_\Gamma(k)\cap\Lambda}
 \left|\Phi^{-1}((\Gamma,\ell))\right| \,. 
\end{equation}
This is a summation over lattice points in the polytope $\Metr_\Gamma(k)$. 
As the $k_i$'s go to infinity, the sum \eqref{sumM} after proper 
scaling will produce an integral.

More precisely, note that, aside from the factor $2^{|k|}$, the right-hand side
of \eqref{e006} is homogeneous in $k$ and $\ell$ of degree $s-\frac32 |e(\Gamma)|$.
Therefore, if the $k_i$'s go to infinity in such a way that
$$
k_i \sim t\cdot\xi_i\,, \quad t\to \infty\,,
$$
the sum \eqref{sumM} becomes the following integral
\begin{multline}\label{eee}
\left|\Map_{g,\Gamma}(k)\right|\sim \\
 \frac{t^{|e(\Gamma)|/2}}{|\Aut(\Gamma)|} 
\frac{2^{|k|-3|e(\Gamma)|/2+1}} {\pi^{|e(\Gamma)|/2}}
\prod_1^s  \xi_i \int_{\Metr_\Gamma(\xi)} d\ell\, 
\prod_{e\in e(\Gamma)} (\ell_{1,e}+\ell_{2,e})^{-3/2} \,,
\end{multline}
where the normalization of Lebesgue measure on 
the polytope $\Metr_\Gamma(\xi)$ is explained in the next
subsection.

The validity of the replacing sums by integrals is justified by
the dominated convergence theorem, 
Lemma \ref{estCat}, and the convergence of the 
following integral
$$
\iint_{\substack{\!\!x,y\ge 0\\\!\!x+y\le c}} \,\,\frac{dx\, dy}{(x+y)^{3/2}} =
\int_0^c \frac{du}{\sqrt{u}}\,, \quad u=x+y \,. 
$$ 

\subsubsection{} 

The Lebesgue measure the right-hand side of \eqref{eee} is normalized
as follows. 

Let $A$ be an open subset of $\Metr_\Gamma(\xi)$. The polytope
$\Metr_\Gamma(\xi)$ is the polytope $\Metr_\Gamma(k)$ scaled by a 
factor of $t^{-1}$ and so $tA \subset \Metr_\Gamma(k)$. 
The number of integer points, that is, the points of the
standard lattice $\Z^{2|e(\Gamma)|}$ in $\Metr_\Gamma(k)$ grows like
$t^{\dim \Metr_\Gamma(k)}$, where the dimension is given by \eqref{dimM}.
We normalize the Lebesgue measure on $\Metr_\Gamma(\xi)$ by 
$$
\int_A  d\ell = \lim_{t\to\infty} t^{s-2|e(\Gamma)|} \, \left|t A\cap \Z^{2|e(\Gamma)|} \right|\,.
$$

The summation in \eqref{sumM} is not over all integer points but over points
in the sublattice $\Lambda\subset \Z^{2|e(\Gamma)|}$ of index $2^{|e(\Gamma)|}$.
Observe that when we intersect both lattices with the affine span of $\Metr_\Gamma(k)$
the index drops to $2^{|e(\Gamma)|-1}$ because one of the parity conditions
becomes redundant once the total perimeter is fixed. Hence
$$
\lim_{t\to\infty} t^{s-2|e(\Gamma)|} \, \left|t A\cap \Lambda \right| =
2^{1-|e(\Gamma)|}\,\int_A  d\ell\,.
$$
This is reflected in
the fact that the exponent of 2 in \eqref{e006} and \eqref{eee} differ by
$|e(\Gamma)|-1$. 

\subsubsection{} 

Consider the sum
$$
\left|\Map_g(k_1,\dots,k_s)\right| =\sum_{\Gamma\in \Gamma^{\ge 3}_{g,s} }
\left|\Map_{g,\Gamma}(k)\right| \,.
$$
Observe, that  some of the summands are asymptotically negligible. Indeed,
it is clear from \eqref{eee} that the asymptotics is determined by those
$\Gamma$ that have the maximal number of edges. Equivalently, 
by invariance of the Euler characteristic, they 
must have the maximal number of vertices. From \eqref{e007} it follows
that this happens if and only if all vertices of $\Gamma$ are
trivalent. 

Denote by $\Gamma^{3}_{g,s}$ the subset of $\Gamma^{\ge3}_{g,s}$
formed by trivalent graphs.  
Remark that every $\Gamma\in\Gamma^{3}_{g,s}$ has $6g-6+3s$ edges.

We have established the following result
\begin{prop}
\begin{multline}\label{e008}
\frac{\map_g(\xi_1,\dots,\xi_s)}{\xi_1\cdots\xi_s}=\\
\frac{2} {(8\pi)^{3g-3+3s/2}}
\sum_{\Gamma\in \Gamma^{3}_{g,s}}
\frac{1}{|\Aut(\Gamma)|} 
\int_{\Metr_\Gamma(\xi)} d\ell\, 
\prod_{e\in e(\Gamma)} (\ell_{1,e}+\ell_{2,e})^{-3/2} \,,
\end{multline}
where $\ell_{1,e}$ and $\ell_{2,e}$ are the lengths of the two sides
of the edge $e\in e(\Gamma)$ in the metric $\ell$.
\end{prop}

\subsubsection{} 

Using the integral
\begin{equation}\label{ab}
\frac{1}{\sqrt{\pi}}\int_0^\infty \int_0^\infty 
\frac{e^{-ax-by}}{(x+y)^{3/2}} \, dx\, dy = 
\frac{2}{\sqrt{a}+\sqrt{b}}\,, \quad \Re a,\Re b >0 \,.
\end{equation}
we can compute the Laplace transform of \eqref{e008} in a compact
form. Take some $z_1,\dots,z_s$ such that $\Re z_i >0$. 

We have 
$$
\int_{\R_{\ge 0}^s} e^{-(z,\xi)}\,\map_g(\xi) \,
\frac{d\xi}{\xi} =
\sum_{\Gamma} \int_{\R_{\ge 0}^s} {d\xi} \int_{\Metr_\Gamma(\xi)} d\ell\,(\dots)
$$
and each summand in the right-hand side is an integral over all possible metrics $\ell$
on $\partial\Gamma$, that is,
just an integral over $\R_{\ge 0}^{2|e(\Gamma)|}$. It factors into a product
of integrals of the form \eqref{ab} over the edges of $\Gamma$.

Thus, we obtain the following 
\begin{thm}\label{t3} The Laplace transform of the function
$\map_g(\xi)$ equals
\begin{multline}\label{ep1}
\int_{\R_{\ge 0}^s} e^{-(z,\xi)}\,\map_g(\xi)\,
\frac{d\xi}{\xi} =
2 \sum_{\Gamma\in \Gamma^{3}_{g,s}}
\frac{1}{|\Aut(\Gamma)|}  
\prod_{e\in e(\Gamma)} 
\frac{2^{-1/2}}{\sqrt{z_{1,e}}+
\sqrt{z_{2,e}}}\,.
\end{multline}
Here $\Gamma^{3}_{g,s}$ is the set of 3-valent ribbon graphs
of genus $g$ with $s$ cells numbered by $1,2,\dots,s$,
$e(\Gamma)$ is the set of edges
of $\Gamma$, and $z_{1,e}$ and $z_{2,e}$ are the two $z_i$'s which correspond
to the two sides of an edge $e\in e(\Gamma)$.
\end{thm}

\subsubsection{}\label{Brown}

The right-hand side is, up to the presence of square roots and
difference in the exponent of $2$, identical
to the right-hand side of the main formula in \cite{K}.
This relation is not accidental. In fact, our counting problem
is very directly related to Kontsevich's combinatorial
description of the intersection numbers on the moduli
spaces. This connection is as follows. 

Consider the following function
\begin{equation}\label{Mh}
\Mh_g(z_1,\dots,z_2) = \sum_{k} e^{-(z,k)} \, \frac{\Map_g(k_1,\dots,k_s)}
{2^{|k|} \prod k_i}\,.
\end{equation}
By definition of $\map_g$, we have
$$
N^{-3g+3-3s/2} \, \Mh_g\left(\frac zN\right) \to \frac 12 
\int_{\R_{\ge 0}^s} e^{-(z,\xi)}\,\map_g(\xi)\,
\frac{d\xi}{\xi}\,,
$$
where the factor $\frac12$ comes from the fact that the summation is
\eqref{Mh} is in fact over all $k$ such that $|k|$ is even which
is an index $2$ sublattice in $\Z^s$. 

It is clear from our discussion that \eqref{Mh} is a sum over
ribbon graphs with vertices of valence $\ge 3$. 
The contribution of each graph $\Gamma$ is
the reciprocal of $|\Aut(\Gamma)|$ times the product
of the following contributions of the
edges $e$ of $\Gamma$. Let $z$ and $w$ be the $z_i$'s 
corresponding to the two sides of $w$, then the contribution $C$ 
of the edge $e$ is
$$
C(z,w)=\sum_{p,q} e^{-zp-qw} \, 2^{-p-q} \, c_{p,q}\,,
$$   
where $c_{p,q}$ is the number of ribbon graphs like the one
shown in Figure \ref{fig8} with length of the upper and
lower boundary being equal to $p$ and $q$ respectively.

Let us forget for a moment that $c_{p,q}$ is just a Catalan
number. Let us think of the Figure \ref{fig8} as of an alley with
trees growing on both sides. The total perimeter of trees
on the two sides is $p$ and $q$.  Let $r$ be the
length of the alley itself, for example, $r=4$ in Figure \ref{fig8}.
Let $t_{p,r}$ denote the number of ways to plant trees of total
perimeter $p$ along an alley of length $r$ so that there
is no tree at the very end of the alley. Clearly
$$
C(z,w)=\sum_r \left(\sum_p e^{-z p}\, 2^{-p} \, t_{p,r} \right) \,
\left(\sum_p e^{-z q}\, 2^{-q} \, t_{q,r} \right) \,. 
$$

\begin{figure}[!hbt]
\centering
\scalebox{.7}{\includegraphics{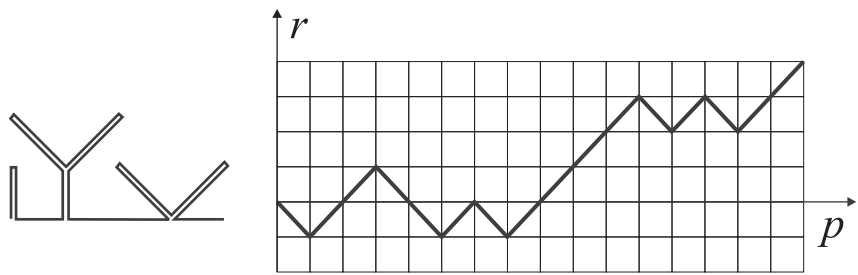}}
\caption{Bijection between planting trees and random walk}
\label{rw}
\end{figure}

It is well known that $t_{p,r}$ also count all trajectories
of a random walk which starting from zero first reach $r$ in
$p$ steps, see Figure \ref{rw}. 
The bijection is very simple: we start a new branch whenever
we go down and finish an existing branch or go to the
next tree whenever we go up. 

Now consider the asymptotics of $C(z/N,w/N)$ as $N\to\infty$.
If $p$ is scaled by $N$ and $r$ by $\sqrt{N}$ then $2^{-p} \, t_{p,r}$
becomes the probability for the standard
Brownian motion to first reach $r$ in time $p$, which
is well known to have the density, see e.g.\ Section V.3.2 in \cite{PR}, 
$$
\frac{r}{\sqrt{2\pi p^3}} \, e^{-r^2/2p} \, dp
$$
with the Laplace transform
\begin{equation}\label{laplt}
\int_0^\infty e^{-z p} \,  \frac{r}{\sqrt{2\pi p^3}} \, e^{-r^2/2p} \, dp  = 
e^{-r\sqrt{2z}} \,.
\end{equation}
Hence
$$
N^{-1/2} \, C(z/N,w/N) \to  
\int_0^\infty e^{-r\left(\sqrt{2z}+\sqrt{2w}\right)} \, dr =
\frac{2^{-1/2}}{\sqrt z + \sqrt w} \,.
$$ 

Kontsevich's combinatorial model for intersection numbers on
the moduli spaces of curves leads to counting alleys with no
trees at all, in which case the length $r$ of the alley is
simultaneously the length of its both boundaries. In our case,
things are dressed up with trees but as \eqref{laplt}
shows this amounts to just replacing
Laplace transform variables by their square roots. 

\subsection{Example: 1-cell maps of genus $\ge 1$}\label{s2ex}

\subsubsection{} 

Consider the case $s=1$, $g=1$. In this case,
the set $\Gamma^{3}_{g,s}$ consists of one element which
is displayed in Figures \ref{fig4} and \ref{fig5}. The automorphism group of this
graph is the cyclic group of order 6 which is clearly seen
in the left half of Figure \ref{fig4}. Also, there is only one $z$ which
corresponds to both sides of every edge. Therefore,
$$
\int_0^\infty e^{-z \xi}\, \map_1(\xi)\, \frac{d\xi}{\xi}
=\frac16\frac1{2^{7/2}}\frac1{z^{3/2}}\,,
$$
which implies that  
$$
\map_1(\xi)=\frac1{12\sqrt{\pi}} \left(\frac\xi2\right)^{3/2}\,.
$$

\subsubsection{}

In general, for $s=1$ and any $g$ we have
$$
\map_g(\xi)\propto \xi^{3g-3/2}\,.
$$
The constant can be fixed using the following exact result of  Harer and Don Zagier \cite{HZ} 
\begin{equation}\label{HDZ}
\left|\Map_g(2k)\right|=\frac{(2k)!}{(k+1)!\,(k-2g)!}\,
\left[x^{2g}\right] \, \left(\frac{x/2}{\tanh x/2}\right)^{k+1}\,,
\end{equation}
where $\left[x^{2g}\right]$ stands for the coefficient of $x^{2g}$. We have
$$
\frac{x/2}{\tanh x/2} = 1+ \frac{x^2}{12} + \dots \,,
$$
which implies that 
$$
\map_g(\xi)=\frac1{\sqrt{\pi}}\frac1{12^g g!}  \left(\frac\xi2\right)^{3g-3/2}\,.
$$

\subsubsection{}

As an exercise, let us check that this is in agreement with \eqref{ymap} and
\eqref{rho1}. In other words, we have to check the identity 
\begin{equation}\label{e009}
\int_{-\infty}^{\infty} e^{\xi x} \, K(x,x) \, dt
= \frac{1}{2\sqrt{\pi}}\, 
\frac{e^{\xi^3/12}}{\xi^{3/2}} \,.
\end{equation}
where $K(x,x)$ is defined by \eqref{rho1}. 

{}From the differential
equation for the Airy function one obtains 
$$
\frac{d^3}{dx^3} K(x,x)-4x \frac{d}{dx} K(x,x) +2 K(x,x)=0\,.
$$
Therefore, its Laplace transform of $K$ must satisfy a first order ODE
which the right-hand side of \eqref{e009} indeed satisfies. This proves
the equality \eqref{e009} up to a constant factor. The factor is fixed
by the asymptotics $\xi\to+ 0$ which was considered in Section \ref{scat}. 

\subsubsection{}

As another application of the exact formula \eqref{HDZ}, let us prove
that taking term-wise limit in
\begin{equation}\label{HDZ3}
\frac{1}{2^k n^{k/2}} \E{\tr H^k} =  \sum_{g\ge 0} n^{1-2g} \, \frac{|\Map_g(k)|}{2^k}
\end{equation}
is justified. We have
$$
\frac{x}{\tanh x} =  1- 2 \sum_{g=0}^\infty (-1)^g \, \zeta(2g) \, \frac{x^{2g}}{\pi^{2g}}
$$
and since $\zeta(2g) < \zeta(2)=\pi^2/6$ for $g>1$, the coefficient of $x^{2g}$ in the above
series is less or equal than $3^{-g}$ in absolute value for any $g$. Therefore, the 
the coefficients of $(1-x^2/12)^{-k-1}$ dominate the coefficients of $\left(\frac{x/2}{\tanh x/2}\right)^{k+1}$
which  implies that 
\begin{equation}\label{HDZ2}
\left|\Map_g(2k)\right|\le   \frac{(2k)!}{(k+1)!\,(k-2g)!} \, \binom{k+g}{g} \, 
12^{-g} \le \frac{1}{\sqrt \pi}\,\frac{2^{2k}\, k^{3g-3/2}}{g!}\,,
\end{equation}
where in the second inequality we used Lemma \ref{estCat} and the inequality $2g\le k$
which implies that 
$$
\binom{k+g}{g}\le \frac{(\frac 32\,k)^g}{g!}\,.
$$

The inequality \eqref{HDZ2} justifies taking term-wise
asymptotics in \eqref{HDZ3} provided 
$k\propto n^{2/3}$ and also yields that
\begin{equation}\label{HDZ4}
\frac{1}{2^k n^{k/2}} \E{\tr H^k}  \le \frac{2^{3/2}}{\sqrt{\pi}} \, \frac{e^{\xi^3/8}}{\xi^{3/2}}\,,
\quad 
\xi= k\, n^{-2/3} \,.
\end{equation}
Using \eqref{sigma_rho} and \eqref{rhorho} we obtain 
\begin{equation}\label{HZest}
\frac1
{2^{|k|}n^{|k|/2}}\E{\prod_{j=1}^s \tr H^{k_i}}
\le \textup{some function of $\xi_i=k_i\,n^{-2/3}$}\,.
\end{equation}
This estimate justifies taking the term-wise limit in \eqref{e22}. 

\section{Random permutations and coverings}

\subsection{Jucys-Murphy elements}
\subsubsection{Definition}

Consider
the following elements $X_1,X_2,\dots$ of the
 group algebra of the symmetric group $S(n)$
\begin{alignat*}{3}
X_1&=(12)+&&(13)+&&(14)+(15)+\dots \,,\\
X_2&=&&(23)+&&(24)+(25)+\dots \,,\\
X_3&=&&&&(34)+(35)+\dots \,,
\end{alignat*}
and so on. These elements are called the Jucys-Murphy elements, or
JM elements for short. For a modern introduction to their properties the reader is
referred to \cite{OV}. See also for example \cite{B,KO,O1,O2} for
various applications of these elements. 

These elements are truly remarkable. Most importantly, they
commute and generate
a maximal commutative subalgebra in the group algebra of $S(n)$
which is exactly the algebra of elements acting diagonally 
 the Young basis of irreducible representations of $S(n)$. 
Since this fact is central to what
follows, we will review it briefly.

\subsubsection{Eigenvalues}

Let $\la$ be a partition of $n$ and consider the corresponding
representation of the $S(n)$. The eigenvalues  of the 
self-adjoint element $X_1$ in the representation $\la$
correspond to the corners of the diagram $\la$ as follows.

\begin{figure}[!hbt]
\centering
\scalebox{.7}{\includegraphics{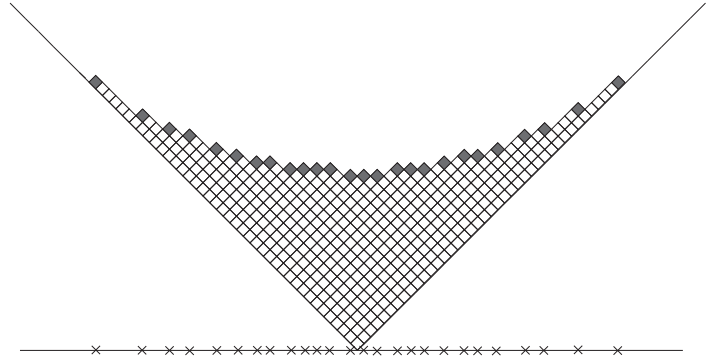}}
\caption{Eigenvalues of the $X_1$}
\label{fig10}
\end{figure}

Let a square $(i,\la_i)\in\la$ be a corner of the diagram $\la$ 
which means that $\la_i>\la_{i+1}$. Then $\la_i-i$ is an eigenvalue
of $X_1$. Recall that the difference between the column number
and the row number of a square $\sq\in\la$ is called the
\emph{content} of $\sq$. That is, the eigenvalues of $X_1$ are
precisely the contents of the corner squares of $\la$. 
If one takes Figure \ref{fig1} and adds the
eigenvalues of the $X_1$ one obtains Figure \ref{fig10}. 

\subsubsection{Eigenspaces}
Now consider the eigenspaces of $X_1$. The subgroup 
$$
S(n)\supset S_1(n) \cong S(n-1)\,,
$$
of permutations which fix $1\in\{1,\dots,n\}$ 
commutes with $X_1$ and thus
preserves the eigenspaces of $X_1$. 

In fact, the $\la_i-i$ eigenspace of $X_1$ is an irreducible module
over $S_1(n)\cong S(n-1)$. Moreover, as an $S(n-1)$-module it 
corresponds to the diagram
$$
\la-\sq_i = (\la_1,\la_2,\dots,\la_i-1,\dots)
$$
obtained from $\la$ by removing the square $\sq_i=(i,\la_i)$. 
In particular, the 
multiplicity of the  eigenvalue $\la_i-i$ equals the dimension
$\dim(\la-\sq_i)$. 

\subsubsection{Action in the regular representation}

Consider the action of $X_1$ in the regular representation, that is,
the representation of $S(n)$ by multiplication on the group algebra
$\C S(n)$. Since the multiplicity of every representation $\la$ in
$\C S(n)$ equals its dimension we find that
\begin{align}
\frac1{n!}\,\tr X_1^k &= \frac1{n!}\,\sum_{|\la|=n} \dim\la \sum_{i} \dim(\la-\sq_i) \,
(\la_i-i)^k \notag \\
&= \sum_{|\la|=n} \Pl_n(\la) \sum_{i} \dn_i(\la)  \,
(\la_i-i)^k \label{dim2} \,,
\end{align}
where the trace is taken in the regular representation,
we agree that $\dim(\la-\sq_i)=0$ if the
square $\sq_i$ is not a corner of $\la$, and we set, by
definition, 
\begin{equation}\label{dn}
\dn_i(\la) =  \frac{\dim(\la-\sq_i)}{\dim\la} \,. 
\end{equation}

The purpose of introducing the ratio \eqref{dn} 
is that it is much simpler than both
its numerator and denominator. Indeed, from the
formula 
$$
\dim\la = |\la|! \, \left. 
\prod_{i<j\le \ell(\la)} (\la_i-\la_j+j-i) \right/
\prod_{i\le\ell(\la)} (\la_i+\ell(\la)-i)\,!\,,
$$
where $\ell(\la)$ is the \emph{length} of the partition $\la$, that is,
the number of nonzero parts in $\la$, it follows that
\begin{equation}\label{dnf}
\dn_i(\la) = \frac{\la_i+\ell(\la)-i}{|\la|} \, \prod_{j\le \ell(\la),j\ne i}
\left(1-\frac{1}{\la_i-\la_j+j-i}\right) \,.
\end{equation}
In the next subsection we will investigate the behavior
of $\dn_i(\la)$ for a $\Pl_n$-typical $\la$ as $n\to\infty$. 

\subsection{Growth and decay of partitions} 

\subsubsection{Rates of growth and decay} 

It is clear that
\begin{equation}\label{res}
\dim \la = \sum_i \dim (\la-\sq_i)\,,
\end{equation}
and hence $\dn_i(\la)$ is naturally a probability
measure on the corners of the diagram $\la$, that is,
$$
\sum_i \dn_i(\la)=1 \,,\quad  \dn_i(\la)\ge 0 \,. 
$$
One can construct a Markov process on the set of all
partitions with transition probabilities 
$$
\Prob\left\{\la \mapsto \la-\sq_i\right\} =\dn_i(\la) \,. 
$$
The representation-theoretic meaning of this decay process
is the branching of representations of symmetric group
under restriction onto a smaller symmetric group. 

Conversely, induction of representations gives a natural
Markov growth process for partitions with transition 
probabilities 
$$
\Prob\left\{\la \mapsto \la+\sq_i\right\} =\up_i(\la) \,, 
$$
where
$$
\up_i(\la) = \frac1{|\la|+1} \, \frac{\dim(\la+\sq_i)}{\dim\la} \,.
$$
We recall that 
the induction rule for representations of symmetric groups implies
that
\begin{equation}\label{ind}
\dim \la = \frac{1}{|\la|+1}\sum_i \dim (\la+\sq_i)\,,
\end{equation}
and hence $\up_i(\la)$ is indeed a probability measure on $i$ for any $\la$. 
Geometrically those $i$ for which $\up_i(\la)\ne 0$ correspond to
places where one can add a square to $\la$, that is, to inner
corners of the diagram $\la$. 

The probabilities $\up_i(\la)$ and $\dn_i(\la)$ are usually called
the transition and cotransition probabilities, see for example \cite{Ke3,Ke5}.
We will call them the rates of growth and decay, respectively.

\subsubsection{Asymptotics of growth/decay rates}

The equations \eqref{res} and \eqref{ind} lead to the following
conclusion: the decay
and growth processes take the Plancherel measure on partition
of $n$ to the Plancherel measure on partitions of $n-1$ and $n+1$,
respectively. 

We are interested in the asymptotics of $\dn_i(\la)$ for fixed $i$
and $\la$ being a Plancherel typical partition of $n$, $n\to\infty$. 
First, let us obtain this asymptotics heuristically. 

Recall that for a Plancherel typical
partition $\la$ of $n$ we have $\la_i \sim 2 \sqrt{n}$, $n\to\infty$. 
This means that after $k$ iterations of the decay process, the 
length of the $i$-row will be $\sim 2 \sqrt{n-k}$. Hence, 
the probability $\dn_i(\la)$ to remove a square from the $i$-th row 
of $\la$ should be
$$
\dn_i(\la) \approx
\frac{2 \sqrt{n} - 2 \sqrt{n-k}}{k} \approx  \frac{1}{\sqrt{n}} \,, \quad n\to
\infty \,.
$$ 
Now let us give a rigorous derivation of this asymptotics

\begin{prop}\label{dec}
With respect to the Plancherel measure $\Pl_n$ on partitions $\la$ of $n$,
$$
\sqrt{n}\, \dn_i(\la) \to 1 \,, \quad n\to\infty\,,
$$
in probability for any fixed $i=1,2,\dots$
\end{prop}

Let us begin with $i=1$. First we show that 
\begin{equation}\label{est1}
\Pl_n\left(\left\{\la, \sqrt{n}\, \dn_i(\la) \le 1+\eps\right\}\right)
 \to 1 \,,
\quad n\to\infty
\end{equation}
for any $i=1,2,\dots$ and any $\eps>0$. Recall that $\la_1\sim
2\sqrt{n}$ and similarly $\ell(\la)\sim 2\sqrt{n}$ for a
Plancherel typical $\la$. Hence,  
$$
\la_1+\ell(\la)\sim 4 \sqrt{n}\,,
$$
and from \eqref{dnf} we obtain 
\begin{equation}\label{dn2}
\dn_1(\la) \sim \frac{4}{\sqrt n} \, \prod_{1<j\le \ell(\la)}
\left(1-\frac{1}{\la_1-\la_j+j-1}\right) \,.
\end{equation}
Note that each factor in the \eqref{dn2} is $<1$. 

The existence of the limit shape $\Omega$ of a typical $\la$ implies
that for any $[a,b]\subset [-2,2]$ the number of $\la_i$'s such that
$\la_i-i\in\sqrt{n}\,[a,b]$ is asymptotic to 
$$
\frac1{\sqrt n}\,
\left|\left\{i,a\le \frac{\la_i-i}{\sqrt n} \le b\right\}\right| \sim 
\frac{\Omega(a)-a}{2}- \frac{\Omega(b)-b}{2} \,, \quad n\to\infty\,.
$$ 
The product over all such $i$ in \eqref{dn2} can be estimated  from
above by
$$
\exp\left(-\frac{1}{2-b}\left(\frac{\Omega(a)-a}{2}- \frac{\Omega(b)-b}{2}
\right)\right) \,.
$$ 
Hence, for any $\eps>0$ we have 
$$
\delta_1(\la) <  \frac{4+\eps}{\sqrt n} \, \exp\left(-\int_{-2}^2 \frac1{2-x}
\,\frac{1-\Omega'(x)}{2}\, dx\right)\,,
$$
for a Plancherel typical $\la$ as $n\to\infty$.  

This integral can be evaluated explicitly  and one finds that 
$$
 \exp\left(-\int_{-2}^2 \frac1{2-x}
\,\frac{1-\Omega'(x)}{2}\, dx\right) = \frac 14 \,.
$$
Indeed, since 
$$
\frac{1-\Omega'(x)}{2} = \frac{1}{\pi}\, \arccos\left(\frac x2 \right) 
$$
one has  to show that
\begin{equation}\label{int1}
\int_{-1}^{1} \frac{\arccos x}{1-x}\, dx = 2\pi\ln 2 \,.
\end{equation}
Changing variables and integrating by parts we obtain
\begin{multline*}
\int_{-1}^{1} \frac{\arccos x}{1-x}\, dx =
\int_0^\pi \frac{t\, \sin t}{1-\cos t} \, dt 
=\pi \ln 2 - \int_0^\pi \ln(1-\cos t) \, dt \,.
\end{multline*}
Using the Fourier expansion 
$$
 \ln(1-\cos t) = \ln((1-e^{it})(1-e^{-it})/2) =
-\ln 2 - 2 \sum_{k=1}^\infty \frac{\cos kt}{k} 
$$
we obtain 
$$
\int_0^\pi \ln(1-\cos t) \, dt = -\pi\ln 2\,,
$$
which establishes \eqref{int1} and \eqref{est1}. 

Observe that, by definition of $\up_1(\la)$, we have  
$$
\up_1(\la)= \frac1{n+1} \, \frac1{\dn_1(\la+\sq_1)}\,, \quad |\la|=n \,.
$$ 
Since $\la+\sq_1$ has the same limit shape $\Omega$, we obtain from
\eqref{est1} that 
\begin{equation}\label{est2}
\Pl_n\left(\left\{\la, \sqrt{n}\, \up_1(\la) \ge  1-\eps\right\}\right)
 \to 1 \,,
\quad n\to\infty
\end{equation}
Now observe that
$$
\sum_{|\la|=n}  \up_1(\la) \, \Pl_n(\la) =
\sum_{|\la|=n}\frac{ \dim\la \, \dim(\la+\sq_1)}{(n+1)!} =
\sum_{|\la|=n+1}  \dn_1(\la) \, \Pl_{n+1}(\la) \,.
$$ 
This, together with \eqref{est1} and \eqref{est2} implies the
existence of both limits 
$$
\sqrt{n}\, \dn_1(\la), \sqrt{n}\, \up_1(\la) \to 1 \,, \quad n\to\infty\,,
$$
in probability.  

Now consider the case $i=2$. First, show that $\la_1-\la_2\to\infty$
for typical $\la$ as $n\to\infty$. Indeed, the formula \eqref{dn2}
can be rewritten as
\begin{equation}\label{dn3}
\dn_1(\la) \sim \frac{4}{\sqrt n} \,  \left(1-\frac{1}{\la_1-\la_2+1}\right)
\, \prod_{2<j\le \ell(\la)}
\left(1-\frac{1}{\la_1-\la_j+j-1}\right) \,.
\end{equation}
It is clear that our analysis of \eqref{dn2} really applies to the
last factor in \eqref{dn3} which means that 
$$
\delta_1(\la) <  \frac{1+\eps}{\sqrt n} \left(1-\frac{1}{\la_1-\la_2+1}\right)
$$
for typical $\la$. Since $\sqrt{n}\, \dn_1(\la) \to 1$ it follows that 
\begin{equation}\label{est3}
\Pl_n\left(\left\{\la, \la_1-\la_2 > \const \right\}\right)
 \to 1 \,,
\quad n\to\infty
\end{equation}
for any constant. 

Therefore we can neglect $(\la_1-\la_2+1)^{-1}$ and write 
\begin{equation*}
\dn_2(\la) \sim \frac{4}{\sqrt n} \, \prod_{2<j\le \ell(\la)}
\left(1-\frac{1}{\la_2-\la_j+j-2}\right) \,.
\end{equation*}
We can apply to this formula exactly the same argument that
we applied to \eqref{dn2} to show that
$$
\sqrt{n}\, \dn_2(\la) \to 1 \,, \quad n\to\infty\,,
$$
in probability. 

An identical argument proves that $\sqrt{n}\, \dn_i(\la) \to 1$
for any fixed $i$. 

\subsection{Plancherel averages and coverings} 

\subsubsection{}
The formula \eqref{dim2} can be rewritten as 
\begin{equation}\label{tr1}
\frac1{2^k \, n^{(k-1)/2} \, n!}\,\tr X_1^k = \sum_{|\la|=n} \Pl_n(\la) \sum_{i} 
\sqrt{n}\,\dn_i(\la)  \,
\left(\frac{\la_i-i}{2\sqrt n}\right)^k \,.
\end{equation} 
Consider the asymptotics of \eqref{tr1} 
as $n,k\to\infty$ in such a way that $k n^{-1/3} \to \xi$ for some
fixed $\xi$. 

The ratio
$
\dfrac{\la_i-i}{2\sqrt n}
$
is maximal (and $\approx \pm 1$) near the edges of the limit shape $\Omega$, that is,
for $i=1,2,\dots$ and also for $i=\ell(\la),\ell(\la)-1,\dots$. We proved
that $\sqrt{n}\,\dn_i(\la)\to 1$ for $i=1,2,\dots$. The condition $k\propto n^{1/3}$
implies that
$$
\left(\frac{\la_i-i}{2\sqrt n}\right)^k \sim \left(\frac{\la_i}{2\sqrt n}\right)^k \,,
\quad i=1,2,\dots\,.
$$ 

What is happening on the other edge of the limit shape is best 
described using the invariance of the Plancherel measure under
$$
\la\mapsto\la'\,,
$$
where $\la'$ denotes the conjugate partition. We conclude that 
\begin{multline}\label{tr2} 
\frac1{2^k \, n^{(k-1)/2} \, n!}\,\tr X_1^k \sim \\
\sum_{|\la|=n} \Pl_n(\la) 
\left(\sum_{i} \left(\frac{\la_i}{2\sqrt n}\right)^k + 
(-1)^k\, \sum_{i} \left(\frac{\la'_i}{2\sqrt n}\right)^k \right)\,.
\end{multline} 
This is totally analogous to the way maximal and minimal eigenvalues
of a random matrix contribute to the asymptotics of  \eqref{e21}. 

\subsubsection{Joint spectrum of JM elements} 

The description of the spectra of the Jucys-Murphy elements $X_i$'s
can be easily iterated. Recall that the $\la_i-i$ eigenspace of 
$X_1$ is the irreducible module over $S_1(n)\cong S(n-1)$ corresponding
to the partition $\la-\sq_i$. This means that the eigenvalues
of $X_2$ in this eigenspace correspond to the corners of the
diagram $\la-\sq_i$ and are irreducible modules over the 
subgroup
$$
S_2(n)\cong S(n-2)\,,
$$
which  fixes $1$ and $2$. Same applies to  $X_3,X_4,\dots$. 

It follows that the formula \eqref{tr2} can be generalized as follows 
\begin{multline}\label{tr3}
\frac1{2^{|k|}\, n^{(|k|-s)/2}\, n!} \tr \prod_{r=1}^s X_r^{k_r} 
\sim \\ 
 \sum_\la \Pl_n(\la)\, 
\left( \sum_{i_1,\dots,i_s=1}^{\infty}
\, \prod_{r=1}^s \left(\frac{\la_{i_r}}{2 \sqrt n}\right)^{k_r}+
\dots
\right) \,,
\end{multline}
where the dots stand for $2^s-1$ more terms involving the 
$\la'_i$'s. Again, this is totally analogous to the situation with 
\eqref{e21}.

\subsubsection{Modified JM elements}

It will be slightly more convenient to consider the following
modification of JM elements. Fix some $s=1,2,\dots$ and set
\begin{align*}
\tX_i&=X_i-\sum_{k=i+1}^s(i\, k)\\
&=(i,s+1)+(i,s+2)+\dots+(i,n)\,, \quad i=1,\dots,s \,.
\end{align*}
We claim that provided $k_i \propto n^{1/3}$ we have
\begin{equation}\label{e010}
 \frac1{n!} \tr \prod_{i=1}^s X_i^{k_i} \sim 
\frac1{n!}  \tr \tX_1^{k_1} \cdots  \tX_s^{k_s}\,, \quad n\to\infty \,.
\end{equation}
Indeed, replacing any $X_i$ with a given transposition leads to the
loss of $\sqrt n$ in the asymptotics because the eigenvalues of
$X_i$ are of order $\sqrt n$. Since there are 
only $|k|\propto n^{1/3}$ possible $X_i$'s to replace, the difference between
the two sides in \eqref{e010} is asymptotically negligible. 

\subsubsection{Traces and equations in $S(n)$}

In the adjoint representation, we have for any $g\in S(n)$ 
$$
\frac1{n!} \tr g =
\begin{cases}
1\,, & g=1\,, \\
0\,, & g \ne 1 \,.
\end{cases}
$$
Hence
$$
\frac1{n!}  \tr \prod_{i=1}^s \tX_i^{k_i} =  \left|\{\tau\}\right|\,,
$$
where $\{\tau\}$ is the set of solutions 
$$
\tau=(\tau_1,\dots,\tau_{|k|})\,, \quad \tau_i\in \{s+1,\dots,n\}\,,
$$
to the following equation in $S(n)$
\begin{equation}\label{etau}
(1\tau_1)\cdots(1\tau_{k_1})(2\tau_{k_1+1})\cdots (2\tau_{k_1+k_2})
\cdots(s\tau_{|k|})=1 \,.
\end{equation}
The symmetric group $S(n-s)$ acts naturally on the set
of all solutions $\{\tau\}$. It is clear that the 
number of elements in the $S(n-s)$-orbit of $\tau$ 
is equal to
$$
\left|S(n-s)\cdot \tau\right| = (n-s)(n-s-1)\dots(n-s-d(\tau)+1) \,,
$$
where $d(\tau)$ is the cardinality of the set 
$\{\tau_1,\dots,\tau_{|k|}\}\subset\{s+1,\dots,n\}$
$$
d(\tau)=\left|\{\tau_1,\dots,\tau_{|k|}\}\right|\,. 
$$
Because $d(\tau)\le |k|\propto n^{1/3}$ we have
$$
\left|S(n-s)\cdot \tau\right|\sim n^{d(\tau)} \,, \quad n\to\infty\,.
$$
It follows that
\begin{equation}\label{sumtau}
\frac 1{ n!} \tr \prod_{i=1}^s \tX_i^{k_i} \sim \sum_{\{\tau\}/S(n-s)} n^{d(\tau)} \,.
\end{equation}

\subsubsection{Equations in $S(n)$ and ramified coverings}
Now remark that the elements of the orbit set $\{\tau\}/S(n-s)$ 
are in bijection with isomorphism classes of
certain coverings of the sphere. 

The corresponding
coverings are defined as follows. Let $0$ be the base point
on the sphere and let $|k|$ points be chosen on a circle around $0$.
It is convenient to assume that $|k|=26$ and denote these points 
by letters of the English alphabet. Our covering  will have 
simple ramifications over $a,b,\dots,z$. That is, the monodromy along 
a small loop encircling each of this point is transposition of 
sheets. 

\begin{figure}[!hbt]
\centering
\scalebox{.7}{\includegraphics{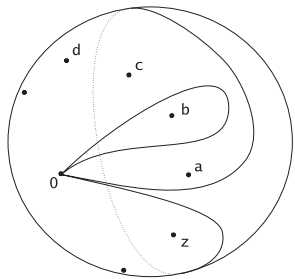}}
\caption{Paths of the monodromy}
\label{fig11}
\end{figure}

In the fiber over
$0$, we pick $s$ sheets, mark them them by $1,\dots, s$,
and call them the {\em special} sheets. We
further require the monodromy around each loop 
around $a,b,\dots$ (see Figure \ref{fig11} where a loop around $b$ is
shown)
to be a transposition of a special sheet with a nonspecial one. 
Another requirement is that 
the first special sheet is permuted by the first $k_1$ loops, the second ---
by the next $k_2$ loops and so on. Finally, we disallow any
unramified sheets. 

The product of all loops, which 
is the big loop in Figure \ref{fig11}, is contractible and so the product of the 
monodromies must be equal to 1.  It is clear, that once we choose
any labeling  of the non-special sheets
in the fiber over 0 by the numbers $\{s+1,\dots,n\}$ 
we get a solution of \eqref{etau} and vice versa. 
Isomorphic coverings differ by a relabeling of the 
the non-special sheets and hence the isomorphism
classes of coverings correspond to $S(n-s)$-orbits. 

We call
the covering satisfying these conditions the 
Jucys-Murphy coverings or {\it JM coverings} for short.
Let $\Su$ be an orientable surface, possibly disconnected. Denote
by $\Cov_\Su(k_1,\dots,k_s)$ the set of JM coverings
$$
\Su \to S^2 \,.
$$ 
It is clear that if a covering corresponds to a solution
$\tau$ of \eqref{etau} then its degree is
$$
\deg \left(\Su \to S^2\right) = s+d(\tau)\,,
$$
and the Euler characteristic of $\Su$ is equal by Riemann-Hurwitz to
\begin{equation}\label{chiSu}
\chi(\Su)=2d(\tau)+2s-|k|\,.
\end{equation}
Therefore, the formula \eqref{sumtau} can be restated as 
\begin{equation}\label{ecov}
\frac1{2^{|k|}\, n^{(|k|-s)/2}\, n!}
\tr \prod_{i=1}^s \tX_i^{k_i} \sim \frac1{2^{|k|}} \sum_{\Su} n^{(\chi(\Su)-s)/2} \, 
\left|\Cov_{\Su}(k_1,\dots,k_s)\right|\,. 
\end{equation}
Here the sum is over all homeomorphism types of orientable
surfaces $\Su$, possibly disconnected. 

As in the case of \eqref{e22}, it is clear that  it is sufficient to
concentrate on connected surfaces only.  If $\Su$ is a connected surface of genus $g$ 
we shall denote the corresponding coverings by $\Cov_g(k_1,\dots,k_s)$. As always,
this set is empty unless $|k|$ is even which we will assume in what follows. 

\section{Counting coverings} 

\subsection{Main result}

\subsubsection{} 

In the present section we will prove the following theorem with
connects JM coverings $\Su\to S^2$ with maps on $\Su$
\begin{thm}\label{t4}
As $k_i\to\infty$, we have
\begin{equation}\label{covmap}
\left|\Cov_g(k_1,\dots,k_s)\right|\sim \left|\Map_g(k_1,\dots,k_s)\right| \,.
\end{equation}
\end{thm}

The proof of this theorem requires some preparations and, in particular,
some understanding of the structure of a JM coverings. Before we start
these preparations, let
us explain how Theorem \ref{t4} implies Theorem \ref{t1}.
Then the rest of the section will be devoted devoted to the proof
of Theorem \ref{t4} and examples.

\subsubsection{Proof of Theorem \ref{t1}}

We know that 
$$
\frac{\left|\Map_g(k_1,\dots,k_s)\right|}{2^{|k|}} \sim  t^{3g-3+3s/2} \, \map_g(\xi)\,, \quad k_i/t\to\xi_i\,,
$$
as $k_i\to\infty$. Hence
if $k_i'\sim n^{1/3} k_i$ then  
$$
\frac{\left|\Cov_g(k_1,\dots,k_s)\right|}{2^{|k|}\, n^{g-1+s/2}} \sim 
\frac{\left|\Map_g(k'_1,\dots,k'_s)\right|}{2^{|k'|} 
\, n^{2g-2+s} }\,. 
$$
I follows that
\begin{equation}\label{covmap2}
\frac{
\left|\Cov_\Su(k_1,\dots,k_s)\right|
}
{2^{|k|} n^{(s-\chi(\Su))/2}}
\sim
\frac{
\left|\Map_\Su(k'_1,\dots,k'_s)\right|
}
{2^{|k'|} n^{s-\chi(\Su)}}
\,,
\end{equation}
provided
\begin{equation}\label{kk}
\frac{k_i}{n^{1/3}}\,,\,\frac{k'_i}{n^{2/3}} \to \xi_i \,, \quad n\to\infty \,.
\end{equation}
It will be clear from the proof of Theorem \ref{t4} that the right
hand side of the following formula \eqref{XH} admits an
estimate of the form \eqref{HZest} and hence we can
apply \eqref{covmap2} termwise to \eqref{e22} and  
\eqref{ecov},\eqref{e010} 
to obtain that 
\begin{equation}\label{XH} 
\frac1{2^{|k|}\, n^{(|k|-s)/2}\, n!}
\tr \prod_{i=1}^s X_i^{k_i} \sim 
\frac1
{2^{|k'|}n^{|k'|/2}}\E{\prod_{j=1}^s \tr H^{k'_i}}
\end{equation}
under the provision \eqref{kk}. Now comparing \eqref{tr3} with the
corresponding result for random matrices finishes the proof. 

Note that the difference in the exponent of $n$ in \eqref{e22}
and \eqref{ecov}
is responsible for the difference in the scaling in 
\eqref{e001} and \eqref{e011}. 

\subsubsection{Proof of Theorem \ref{t1'}}

Our argument follows the argument of Section 5 in \cite{Sosh}.
Consider the following random measures on $\R$
$$
\fX=\sum_i \delta_{x_i} \,, \quad \fY=\sum_i \delta_{y_i}\,,
$$
where
$$ 
x_i=n^{1/3}\left(\frac{\la_i}{2 n^{1/2}}-1\right) \,,
\quad 
y_i=n^{2/3}\left(\frac{E_i}{2 n^{1/2}}-1\right)\,.
$$ 
Define $\E{\fX^{\times s}}$ as the following nonrandom measure
on $\R^s$ 
$$
\E{\fX^{\times s}}
(A_1\times \cdots \times A_s) = \E{\prod_{i=1}^s \fX(A_i)}\,,
$$
and define $\E{\fY^{\times s}}$ similarly. Theorem \ref{t1}
says that the Laplace transforms of $\E{\fX^{\times s}}$
and $\E{\fY^{\times s}}$ have the same limit as $n\to\infty$.
Multiply $\E{\fX^{\times s}}$  and  $\E{\fY^{\times s}}$ by the exponential
of the sum of coordinates, which is equivalent to
shifting Laplace transform variables by $(1,\dots,1)$. This
yields finite measures for which the convergence of Laplace
transforms implies weak convergence. Hence all mixed
moments of the following random variables
\begin{equation}\label{XA}
\fX(A)=|\{y_i\in A\}|\,, \quad \fY(A)=|\{y_i\in A\}|\,, \quad A\subset [c,\infty)^s\,,
\end{equation}
have identical limits, where $c\in\R$ is arbitrary fixed and $A$ varies. 
From this one concludes (cf.\ \cite{Sosh})
that the joint distributions of the random variables \eqref{XA}
are the same in the $n\to\infty$ limit. The theorem follows. 

\subsection{Structure of JM coverings}

\subsubsection{Valence of nonspecial sheets} 

\begin{figure}[!hbt]
\centering
\scalebox{.7}{\includegraphics{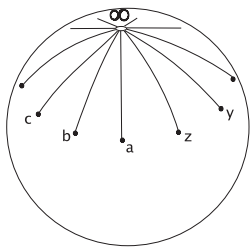}}
\caption{Cuts on the sphere}
\label{fig12}
\end{figure}

 Let us make $|k|$ cuts on
the sphere from the points $a,\dots,z$ to the infinity
as in Figure \ref{fig12}.
This cuts $\Su$ into $(s+d)$ polygons. Let us 
describe the shape of these polygons and how they fit together. 

Given a nonspecial sheet $\sigma$, let its {\em valence} be the
number of points from $a,\dots,z$ such that the monodromy
around that point permutes $\sigma$. Clearly, the
valence of every sheet is $\ge 2$. On the other hand 
\begin{equation}\label{e015}
\sum_{\text{nonspecial $\sigma$}}
 (\val(\sigma)-2)= |k|-2d =2s-\chi(\Su)\,,
\end{equation}
therefore the number of sheets of valence $\ge 3$ is
bounded by $2s-\chi(\Su)$. 

\begin{figure}[!hbt]
\centering
\scalebox{.7}{\includegraphics{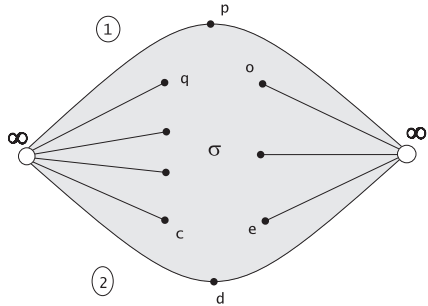}}
\caption{Nonspecial sheet of valence 2}
\label{fig13}
\end{figure}

Suppose $\sigma$ is a 2-valent sheet and suppose that
the monodromy around one ramification point, say, $p$ permutes it with 
the 1st
special sheet and the monodromy around another
ramification point, say,  
$d$ permutes it with 2nd special sheet.  Then the preimages of the
cuts in Figure \ref{fig12} on $\Su$ are drawn in 
Figure \ref{fig13} where the circled numbers 1 and 2 indicate
that the corresponding boundary is attached to the
1st and 2nd special sheet, respectively. 
Note how the angles get halved at the points which
cover the points $p$ and $d$. 

\begin{figure}[!hbt]
\centering
\scalebox{.7}{\includegraphics{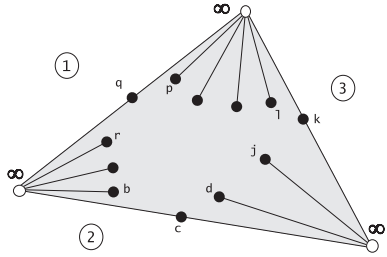}}
\caption{Nonspecial sheet of valence 3}
\label{fig14}
\end{figure}

Similarly, if $\sigma$ is a 3-valent sheet then 
it looks like a triangle (similarly, a sheet of
valence $m$ looks like an $m$-gon). For example
if monodromies around  $q$, $c$, and $k$  permute $\sigma$ with
the the 1st, 2nd, and 3rd special sheet, respectively,
then $\sigma$ looks like Figure \ref{fig14}. 

\subsubsection{Ribbon graph associated to a covering} 

\begin{figure}[!hbt]
\centering
\scalebox{.7}{\includegraphics{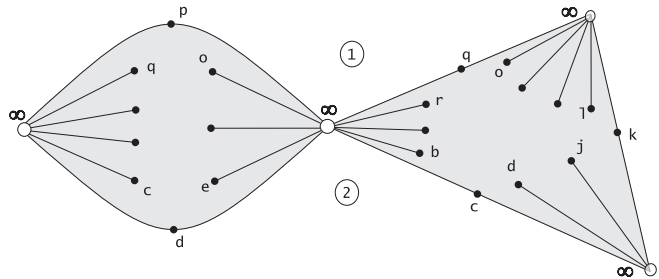}}
\caption{How nonspecial sheets fit together}
\label{fig15}
\end{figure}

The nonspecial sheets naturally glue together at
the points which cover $\infty$ to form a ribbon graph
whose edges are the 2-valent sheets and 
vertices are either the sheets of valence $\ge 3$
or multivalent junctions (like in Figure \ref{fig20})
of 2-valent sheets. See Figure \ref{fig15}
and note how $q$ follows $p$ and $d$ follows $c$ after
passing through $\infty$.  

Observe, in particular how we have the whole alphabet
going once clockwise around each point over $\infty$. 
This reflects the fact that there is no
ramification over $\infty$.

\subsubsection{Special sheets} 

\begin{figure}[!hbt]
\centering
\scalebox{.7}{\includegraphics{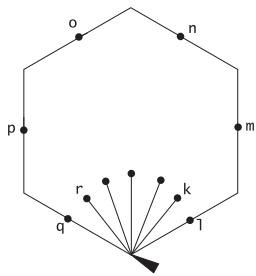}}
\caption{A special sheet with a marked vertex}
\label{fig15a}
\end{figure}

The cells of this ribbon graph correspond to the
special sheets and look as follows. Suppose
$\sigma$ is the $i$-th special sheet. Then the
valence of $\sigma$ is, by construction, equal to $k_i$.
Suppose that $k_i=6$ and the corresponding
ramification points are $\{l,m,m,o,p,q\}$. Then this
special sheet looks like the hexagon in Figure \ref{fig15a}. 
The special sheets come with a natural choice of the 
marked vertex, namely, the initial vertex of their first edge
in alphabetical order. For example, in Figure \ref{fig15a} the 
bottom vertex is the
marked vertex.

\subsubsection{Examples}\label{cox}

\begin{figure}[!hbt]
\centering
\scalebox{.7}{\includegraphics{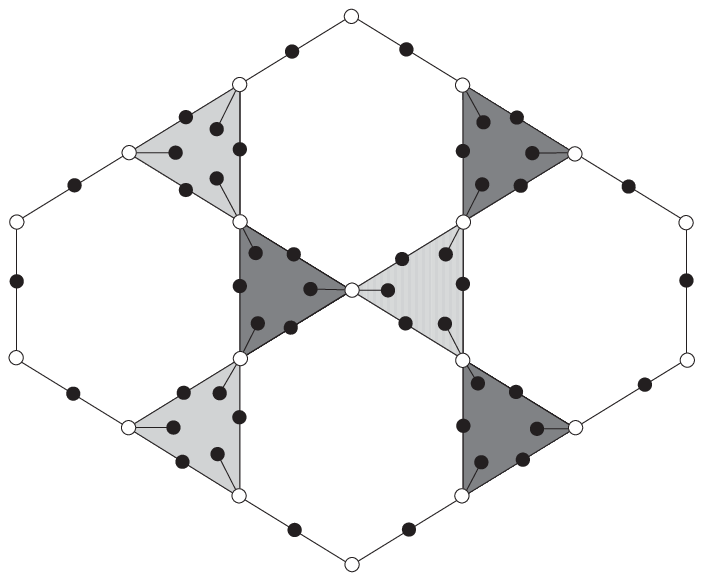}}
\caption{A 3-fold JM covering $T^2\to S^2$}
\label{fig15b}
\end{figure}

Consider the following solution to \eqref{etau}
$$
(12)(13)(12)(13)(12)(13)=1\,,
$$
which is the Coxeter relation in $S(3)$. The corresponding
3-fold covering of the sphere is a torus and the 
3 sheets (one 6-valent special, two 3-valent nonspecial)
fit together on the torus $T^2$ shown in Figure \ref{fig15b}.

\subsection{From coverings to maps} 

\subsubsection{The collapsing mapping $\Psi$} 

We introduce now  the following mapping $\Psi$ from 
JM coverings $\Su\to S^2$  with $s$ special sheets to 
maps on $\Su$ with $s$ boundary components. What $\Psi$ does is it simply
collapses all nonspecial sheets as follows. 

\begin{figure}[!hbt]
\centering
\scalebox{.7}{\includegraphics{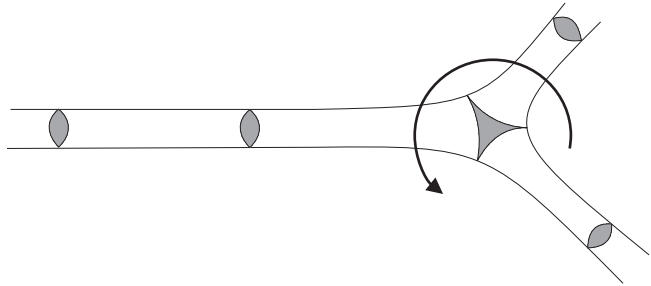}}
\caption{Collapse of Figure \ref{fig15}}
\label{fig16}
\end{figure}

If
a nonspecial sheet $\sigma$ is 2-valent then we plainly
collapse it and glue together the two special sheets which
$\sigma$ separated.  Nonspecial sheets of
 valence $\ge 3$ we shrink to the middle
as shown in Figure \ref{fig16} where the collapse of 
the two nonspecial sheets from Figure \ref{fig15} is
shown (the meaning of the arrow in Figure \ref{fig16}
will be explained below).
Note that collapsing a sheet $\sigma$ of valence $\val(\sigma)\ge 3$
increases the length of each of the $\val(\sigma)$ boundaries involved
by 1. For example, the boundary in Figure \ref{fig16} is 
3 units longer than the boundary in Figure \ref{fig15}. 

The special sheets become the cells of the map, their numbering
is just the numbering of the special sheets by $1,\dots,s$ and
the marked vertices are the marked vertices of the special 
sheets. 

\subsubsection{Example}

Note that collapsing the covering discussed in
Section \ref{cox} and shown in Figure \ref{fig15b}
produces, essentially,
 the map on torus shown in Figures \ref{fig4} and \ref{fig5}. 
More precisely, every edge of this map has length 2 instead
of $1$, so the torus is really glued from a 12-on, not from 
a hexagon.

\begin{figure}[!hbt]
\centering
\scalebox{.4}{\includegraphics{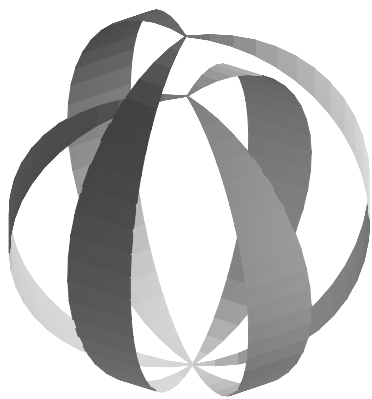}}
\caption{Ribbon graph corresponding to $(12)^6=1$}
\label{(12)6}
\end{figure}

\begin{figure}[!hbt]
\centering
\scalebox{.5}{\includegraphics{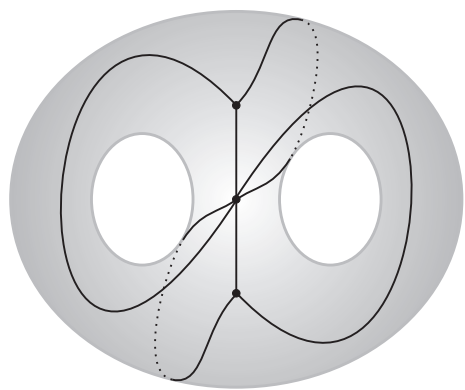}}
\caption{Embedding of the graph form Figure \ref{(12)6}}
\label{genus2}
\end{figure}

As another example, consider the equation
$$
(12)^6 = 1 \,.
$$
Which defines a $2$-fold covering of the sphere of genus $2$
with $1$ special and $1$ nonspecial sheet, both $6$-valent.
This nonspecial sheet look like a hexagon with 3 nonadjacent
vertices glued together and 3 other nonadjacent
vertices also glued together.  When we collapse
this figure to the middle to get the ribbon graph 
shown in Figure \ref{(12)6}. 
Its embedding into the genus 2 surface is shown in 
Figure \ref{genus2}

\subsubsection{Left and right vertices}

We now observe that ribbon graphs associated to JM covering and
the corresponding maps have vertices
of two following fundamentally different types. 
Let $v$ be a vertex of a map.
Suppose we are going around the boundary of the 1st polygon
counterclockwise, then the around the boundary of the 2nd 
polygon counterclockwise and so on. We visit our vertex
$\val(v)$ times from the $\val(v)$ corners which meet
at $v$. We call the vertex $v$ a {\it right} vertex if the
corners are visited in the clockwise order and a 
{\it left} vertex if the corners are visited in the 
counterclockwise order. By definition, we call $v$ right
if $\val(v)\le 2$. 

\begin{figure}[!hbt]
\centering
\scalebox{.7}{\includegraphics{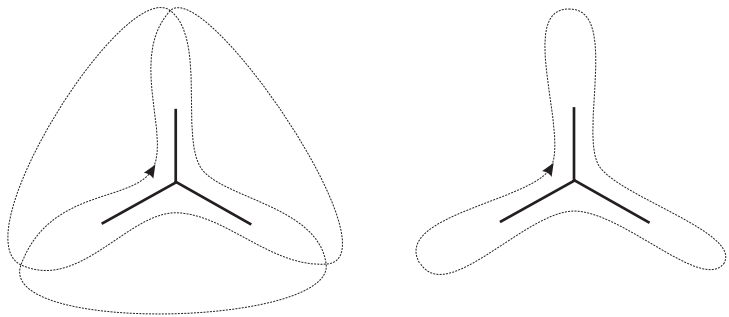}}
\caption{A left vertex and a right vertex}
\label{fig16b}
\end{figure}

Note that if $\val(v)>3$ then
$v$ may be neither left nor right. For an example of this,
look at the surface of genus 2 obtained by
identifying opposite sides of a 10-gon. A left and right
vertex of $\val(v)=3$ are shown in Figure \ref{fig16b} where the 
dashed lines represent the order of going around the 
three corners. 

Suppose $v$ is a vertex of map which came from of
a JM covering. Then $v$ either covers $\infty$ or
$v$ is the middle point of a collapsed $m$-valent
nonspecial sheet where $m\ge 3$. Observe that then
$v$ is a right or left vertex, respectively. Indeed,
if $v$ covers $\infty$ then, since there is no
ramification over $\infty$, the whole alphabet is
circling $v$ once clockwise. Similarly, if $v$ was
a midpoint of a nonspecial sheet then (see Figure \ref{fig14})
the alphabet was going around $v$ once counterclockwise.
This translates into $v$ being a right and left
vertex, respectively. 

For example, the arrow in Figure \ref{fig16}
shows the order of visiting the corners of the
trivalent vertex (which is left). In  another
example, the graph from Figures \ref{(12)6} and
\ref{genus2} has one  6-valent left vertex and two
 3-valent right vertices.

\subsubsection{The image of $\Psi$}
We call a vertex $v$ an {\em interior} vertex 
if all corners which at meet $v$ come from the 
same polygon of the map. 

We will now prove the following

\begin{prop}\label{p5}
The mapping $\Psi$ from JM coverings
to maps is one-to-one.  Its image $\Img\Psi$ consists of all
maps satisfying the two following conditions:
\begin{itemize}
\item every vertex is either left of right,
\item all marked vertices are interior right vertices,
\item the distance between any two left vertices is $\ge 2$. 
\end{itemize}
\end{prop}

Recall that for vertices of valence $\ge 4$ being left
or right is a nontrivial condition and that is was 
shown above that only left or right vertices arise
from JM coverings. 

\medskip 

\begin{figure}[!hbt]
\centering
\scalebox{.7}{\includegraphics{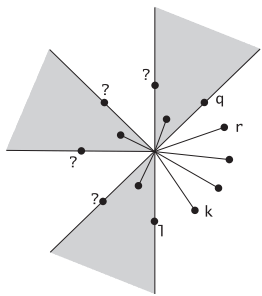}}
\caption{The marked vertex must be an inner vertex}
\label{fig16c}
\end{figure}

\noindent
{\em Proof.}
By construction, all marked vertices come from
some points which cover $\infty$ and, therefore, 
they are right vertices. Let us show that they also
must be interior vertices. This follows from inspection
of Figure \ref{fig16c}. The Figure \ref{fig16c} shows the marked 
vertex (the bottom one) of the special sheet from Figure \ref{fig15a}. 

Since the whole alphabet must go  once around $\infty$
the points marked by question marks in Figure \ref{fig16c}
cannot be points from $\{r,s,\dots,j,k\}$.
On the other hand, the points $\{r,s,\dots,j,k\}$ are
precisely the ramification points which {\em do not}  
lie on the boundary of our special sheet. 
Therefore, all points marked by question marks
do lie on the boundary of our special sheet.
It follows that all corners in Figure \ref{fig16c} come come from one and
the same special sheet. 

An algebraic equivalent of this geometric argument is the
following. Let $\tau$ be a solution of 
\eqref{etau}. Then $(1\tau_1)\cdots(1\tau_{k_1})$ must fix $1$ because
$1$ is clearly fixed by the rest of the product in \eqref{etau}. This
translates into Figure \ref{fig16c}.

We will now show that any map satisfying the above 
two conditions comes from a unique JM covering. This covering can be 
reconstructed as follows. 

Assign symbols $a,b,c,\dots$ consecutively to all edges 
of the polygons of the map starting from 
the marked vertex of the first polygon. 

Now consider some vertex $v$ of our map. If $v$  is  
a right vertex (in particular, if $\val(v)\le 2$)
then the structure of the corresponding
JM covering at $v$ can be reconstructed uniquely
from the fact that $v$ covers $\infty$ and there is
no ramification at $\infty$. (In other words, all letters of
the alphabet have to occur once clockwise around $v$).

\begin{figure}[!hbt]
\centering
\scalebox{.7}{\includegraphics{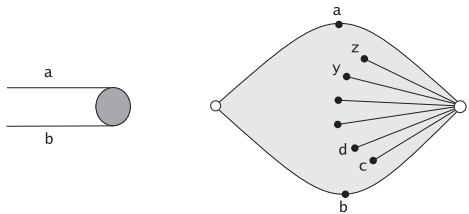}}
\caption{Reconstruction of JM covering for $\val(v)=1$}
\label{fig17}
\end{figure}

\begin{figure}[!hbt]
\centering
\scalebox{.7}{\includegraphics{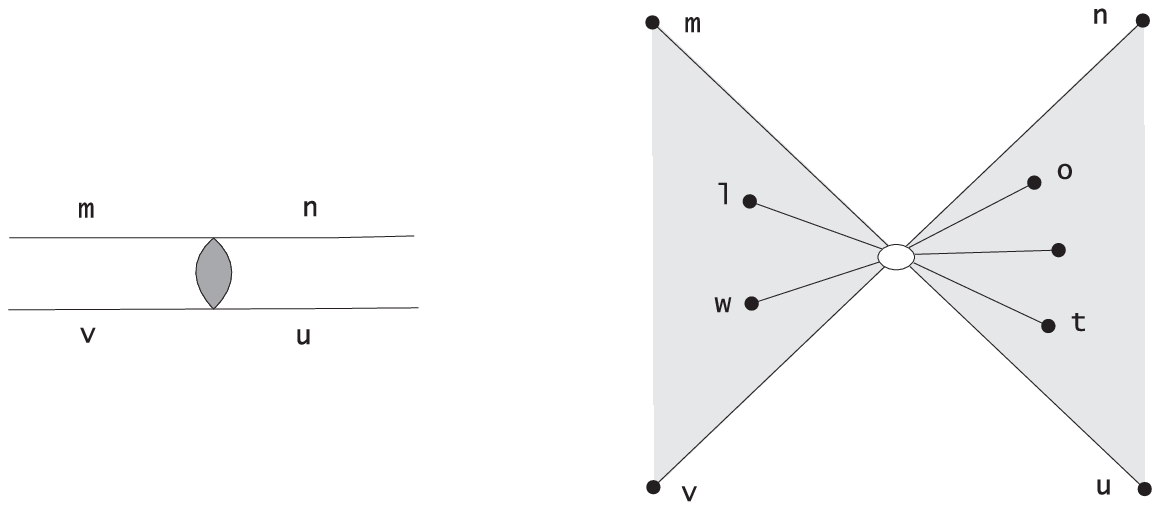}}
\caption{Reconstruction of JM covering for $\val(v)=2$}
\label{fig18}
\end{figure}

\begin{figure}[!hbt]
\centering
\scalebox{.7}{\includegraphics{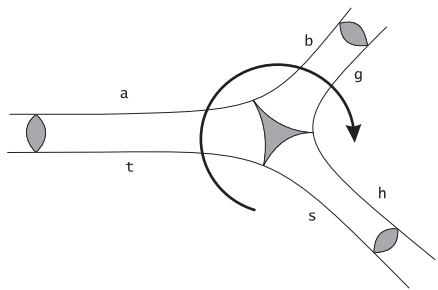}}
\caption{A right vertex of valence $3$}
\label{fig19}
\end{figure}

\begin{figure}[!hbt]
\centering
\scalebox{.7}{\includegraphics{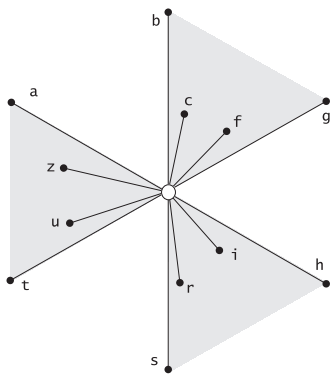}}
\caption{The JM covering corresponding to a right vertex of valence 3}
\label{fig20}
\end{figure}

This reconstruction is shown, respectively,  in Figures \ref{fig17}
for the case $\val(v)=1$, in Figure \ref{fig18} for
the case $\val(v)=2$, and in Figures \ref{fig19} and
\ref{fig20} for $\val(v)=3$. 

When we encounter a left vertex (such as the vertex with the
arrow in Figure \ref{fig16}) then  we insert a nonspecial
sheet of valence $\val(v)$. The result looks like
Figure \ref{fig14}. One should notice that
this operation reduces the number of edges
by $\val(v)$ and we have to relabel the edges if we
want a consecutive alphabetical labeling. 
Also notice that we have to use the third condition
in the statement of the proposition for this
reconstruction. 

This finishes the reconstruction of
nonspecial sheets. As to the special sheets,
let us examine the Figure \ref{fig16c}. 
If the edges of a cell of a map
are labeled by $l,m,\dots,p,q$ and its initial
vertex is an interior right vertex then there
is room to fit in the rest of the alphabet 
as in Figure \ref{fig16c}.  This concludes the proof. 

\subsubsection{Example: $s=1$ and $g=0$} 
Consider the case $s=1$ and $g=0$. The equation \eqref{e015}
implies that in this case
$$
\sum_{\text{nonspecial $\sigma$}}
 (\val(\sigma)-2)= 2s-2+2g=0\,,
$$
and hence there are no nonspecial sheets of valence $>2$.
Therefore the map $\Psi$ is a bijection between the sets 
$\Cov_0(k)$ and $\Map_0(k)$.

The algebraic translation of this geometric fact is the 
following. Let 
\begin{equation}\label{e017}
(1i_1)(1i_2) \dots (1i_k) = 1
\end{equation}
be the solution of \eqref{etau} corresponding to our covering. 
By \eqref{chiSu} the condition $s=1$ , $g=0$ implies that
$$
2 d(\tau) - |k| = 2-2g-2s = 0 
$$
and since every $i_j$ has to appear at least twice, this
is equivalent to saying that there are precisely $k/2$ pairs of 
equal numbers among the numbers $i_1,\dots,i_k$.

The bijection $\Psi$ between $\Cov_0(k)$ and $\Map_0(k)$ and the 
example in Section \ref{scat} now mean that \eqref{e017} is 
satisfied if and only if the the corresponding pairing is noncrossing.
This observation is due to P.~Biane \cite{B2}.

Note  that the noncrossing in \eqref{e017} means that this
equality is a consequence of solely the relations
$$
(1i)^2 = 1\,, \quad i=1,2,\dots\,,
$$
among the generators $(12),(13),(14),\dots$ of the symmetric
group.

\subsection{Counting maps}

\subsubsection{} 

Intoduce the following subsets in $\Map_g(k)$, where as usual, 
we use the abbreviation
$$
k=(k_1,\dots,k_s)\,.
$$
Denote by
$$
\Map^3_g(k)=\Phi^{-1}\left(\Gamma^{3}_{g,s}\right)
$$
the set of those maps which after contraction $\Phi$ have
only trivalent vertices.   Since only trivalent graphs 
$\Gamma$ contribute to \eqref{e008}, we know that 
\begin{equation}\label{e014}
\left|\Map^3_g(k) \right| \sim 
\left|\Map_g(k)\right|\,, 
\quad k_i\to\infty\,.
\end{equation}

By definition, set  
$$
\Map^*_g(k) = \Map^3_g(k) \, \cap\, \Img\Psi\,,
$$
that is, $\Map^*_g(k)$ is 
the subset of $\Map^3_g$ formed by maps
satisfying the conditions of Proposition \ref{p5}.
It is the image under $\Psi$ of JM 
coverings with nonspecial sheets of valence at most $3$.

We will establish the following 

\begin{prop}\label{p6} 
\begin{align}\label{e012}
&\left|\Map^*_g(k) \right|
\sim 2^{-6g+6-6s} \left|\Map^3_g(k)\right|
\\\label{e013} 
&\left|\Img\Psi\cap \Map^3_g
\setminus \Map^*_g\right| = 
o\left(\left|\Map^*_g\right|\right) \,.
\end{align}
\end{prop}

Once Proposition \ref{p6} is established, the Theorem \ref{t4} will
follow. Indeed, since $\Psi$ is one-to-one, 
then because of \eqref{e014} and \eqref{e013}
 it suffices to consider coverings
with only $\le 3$-valent nonspecial sheets. For such a covering, 
the number of 3-valent sheets equals 
$2g-2+2s$. Collapsing a trivalent sheet to its middle
increases the length of the boundary by 3. Therefore, in total,
the boundary of the corresponding map is $6g-6+6s$ longer.
Since this precisely compensates the exponent in \eqref{e012}, we obtain: 
$$
\left|\Cov_g(k) \right| \sim 
\left|\Map_g(k)\right|\,, 
\quad k\to\infty\,.
$$

We also point out that Proposition \ref{p5} gives an upper bound
on $|\Cov_g(k)|$ which results in an analog of the estimate
\eqref{HZest} for the right-hand side of \eqref{XH}. Indeed,
the mapping $\Psi$ is one-to-one and increases the total perimeter
of the boundary by at most $6g-6+6s$. From \eqref{chiSu} we
have $g\le |k|/2$, so the total increase is by at most a
multiple of $|k|$ with implies that the right-hand side 
of \eqref{XH} is again bounded by some function of the $\xi_i$'s. 

\subsubsection{Proof of Proposition \ref{p6}}

In order to examine the difference between the sets $\Map^3_g$, $\Map^*_g$,
and $\Img\Psi$ we need to introduce the following notions.

\begin{figure}[!hbt]
\centering
\scalebox{.7}{\includegraphics{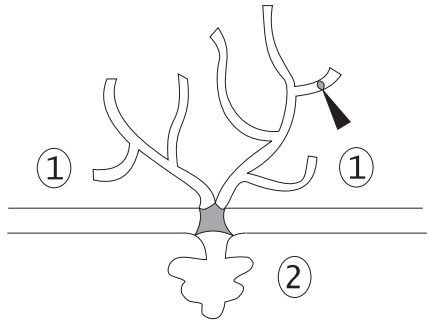}}
\caption{A mouth vertex}
\label{fig22b}
\end{figure}

Let $v_0$ be a
marked vertex of a map with $s>1$ polygons. Suppose
that $v_0$ is an interior vertex. Follow the edges
of the corresponding polygon in the counterclockwise
direction until we reach a vertex $v$ which is not interior. 
By analogy with the flow of a river, we call the vertex $v$
a {\em mouth} vertex, see Figure \ref{fig22b}.
Observe that a mouth vertex is {\em never
right}. 

Also, call a vertex $v$ of a map {\em contractible} if it 
disappears after contraction of all $\le 2$-valent vertices;
otherwise, call it {\em incontractible}. Observe that a
contractible vertex is always right unless it is a mouth
vertex. 

\medskip 

\noindent
{\em Proof.}
First, the condition
that the marked vertices must be right is asymptotically
negligible. Indeed, all but finitely many vertices
are right and the chances to hit one them with
a mark go to 1 as the perimeter goes to infinity.  
Similarly, the third condition in Proposition
\ref{p5} does not affect the asymptotics. 

The possible combinatorial configurations of the
incontractible and mouth vertices of maps in
$\Map^3_g$ are  described by
ribbon graphs $\Gamma\in\Gamma^{3}_g$ together with the
choice of an edge $e_i\in e(\Gamma)$,
$i=1,\dots,s$, on the boundary of any cell of $\Gamma$. The edge $e_i$
is the first edge we reach if we start from the 
marked vertex of the $i$-th cell of the map.
We shall see that,  for any configuration, the proportion of
maps lying in $\Map^*_g$ equals 
$\sim2^{-6g+6-6s}$,   and that the same portion of maps
lies in $\Img\Psi$. This number is,
in fact, a product of factors $2^{-3}$ over 
the $2g-2+s$ vertices $v$ that are not right. Let us examine
such vertices $v$. 

\begin{figure}[!hbt]
\centering
\scalebox{.7}{\includegraphics{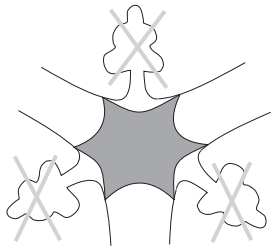}}
\caption{The forbidden trees of a left vertex}
\label{fig23}
\end{figure}

First, suppose $v$ is an incontractible vertex. 
By definition of $\Map^3_g$, it means that $v$ becomes
trivalent after the 3 trees shown in
Figure \ref{fig23} are contracted  onto it. 
We claim that for such a vertex being left is
equivalent to being trivalent. Indeed, suppose
$v$ is left and not trivalent. Then, as we go around 
any nonempty tree in any of the trees shown
in Figure \ref{fig23},  
we go from one corner of $v$ to the next corner
in the clockwise direction. Since $v$ is left, this is
impossible.

It follows from the discussion in Sections \ref{sfor}  or
\ref{Brown} that 
the removal of any given tree comes at the price
of the factor $\frac12$ in the asymptotics. In terms
of the random walk, for example, it means that the first
step of the walk has to go up, which is an event of  probability
$\frac12$. Therefore,
asymptotically $\sim\frac1{2^3}$ of maps are trivalent
at $v$ or, equivalently, $v$ is a left (or trivalent)  vertex for
about $\sim\frac1{2^3}$ of all maps. 

Now suppose that $v$ is a contractible mouth vertex,
such as the one shown in Figure \ref{fig22b}. We may 
assume that $v$ does not coincide with any other
mouth vertex because the chances of 
such a coincidence vanish as the perimeter of the map
goes to infinity.  

With this assumption, $v$ being left is again equivalent
to $v$ being trivalent and both mean that $v$ must look like
the vertex in Figure \ref{fig22}:
\begin{figure}[!hbt]
\centering
\scalebox{.7}{\includegraphics{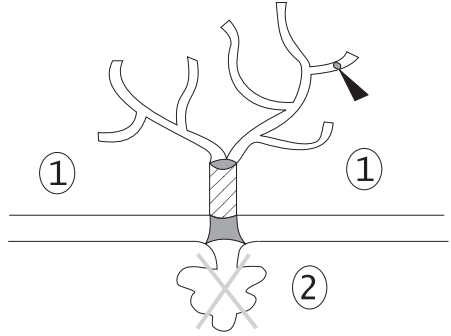}}
\caption{A trivalent contractible mouth vertex}
\label{fig22}
\end{figure}
namely, the tree at the bottom
must be empty and only one branch (shaded in Figure \ref{fig22})
must go up. 

For general maps, multiple branches may go up 
or no branches at all (which happens if the marked vertex
is not interior). Therefore, the insertion of this shaded
branch and chopping down the tree at the bottom takes
arbitrary maps to maps such that $v$ is a trivalent 
mouth vertex and the corresponding marked vertex is interior.
The insertion of the shaded branch increases the perimeter
by $2$. This means that $\sim\frac1{2^2}$ of maps have it.
This times $\frac12$ for the forbidden tree gives us 
the total of $\sim\frac1{2^3}$  of maps belonging to
$\Map^*_g$.

Either way, we get a factor of $2^{-3}$ for any trivalent left vertex. 
The number of such vertices can be easily computed. All of them
become trivalent nonspecial sheets of
the corresponding JM covering. Therefore, by \eqref{e015} there are
$2g-2+2s$ of them. This proves \eqref{e012} and \eqref{e013}
and concludes the proof
of  Proposition \ref{p6} and, hence, of Theorem \ref{t1}.

\subsection{Example}

Note  that the noncrossing in \eqref{e017} meant that this
equality is a consequence of solely the relations
$$
(1i)^2 = 1\,, \quad i=1,2,\dots\,,
$$
among the generators $(12),(13),(14),\dots$ of the symmetric
group. The relations of Coxeter type (which produce coverings
of genus 1)
$$
(1i)(1j)(1i)(1j)(1i)(1j)=1
$$
start playing role in the enumeration of $\Cov_1(k)$. 

Every covering in $\Cov_1(k)$ has
either two 3-valent special sheets or, else, one of valence 4.
Consider the first case because the second makes
no contribution to the asymptotics. Denote by $\Cov^3_1(k)$
the corresponding subset of $\Cov_1(k)$. 

For $\Cov^3_1(k)$, the corresponding
relations are, up to a cyclic shift:
\begin{equation}\label{e018}
(1i) \, w_1 \,(1j)\, w_2 \,(1i)\, w_3\, (1j)\, w_4\, (1i)\, w_5\, (1j)\, w_6 =1
\end{equation}
Here the $w_i$'s are some words in the generators $(12),(13),(14),\dots$
subject to two conditions. First, $(1i)$ and $(1j)$ appear exactly
3 times each in \eqref{e018} and any other generator appears either 0 or 2
times. Second, 
$$
w_1 \, w_4 =w_2 \, w_5 = w_3\, w_6 = 1 \,,
$$
which means that any relation \eqref{e018} is built from 3 relations
from the $g=0$ case. Using the generating function for the
Catalan numbers, one obtains the following 
generating function
$$
\sum_k \left|\Cov^3_1(k) \right|\, z^k = 
=\frac14
\frac{z^2\left(1-\sqrt{1-4z^2}\right)^2}
{(1-4z^2)^{5/2}} \,.
$$
Since
$$
\frac14
\frac{z^2\left(1-\sqrt{1-4z^2}\right)^2}
{(1-4z^2)^{5/2}} \sim 
\frac1{16} \frac1{(1-4z^2)^{5/2}}\,, \quad z^2\to \frac14\,,
$$
we conclude that
$$
\frac{\left|\Cov_1(k) \right|}{2^k} \sim 
\frac{\left|\Cov^3_1(k) \right|}{2^k} \sim
\frac1{16} \frac{(k/2)^{\frac52-1}}{\Gamma(5/2)} =
\frac1{12\sqrt{\pi}} \left(\frac k2\right)^{3/2}\,,
$$
which agrees with computations of Section \ref{s2ex}

\end{document}